\documentclass[twocolumn]{autart}

\usepackage{amsmath,amssymb}
\usepackage{graphics}
\usepackage{epsfig}
\usepackage{graphicx}
\usepackage{epstopdf}
\usepackage[colorlinks=true,breaklinks]{hyperref}
\newtheorem{theorem}{Theorem}[section]

\newtheorem{lemma}[theorem]{Lemma}
\newtheorem{proposition}{Proposition}

\newtheorem{remark}{Remark}
\newtheorem{example}{Example}
\newtheorem{assumption}{Assumption}

\begin{document}
\begin{frontmatter}

\title%[Cycle stabilization with noise]
{Stabilization of cycles for difference equations
with a noisy PF control}

\thanks[footnoteinfo]{The first author was supported by NSERC grant RGPIN-2015-05976,
the second  and the fourth by the Grant FEKT-S-17-4225 of Faculty of 
Electrical Engineering and Communication, Brno University of Technology. The third  
author was supported  by the project International Mobility of Researchers of Brno
University of Technology  CZ.02.2.69/0.0/0.0/16-027/0008371. E. Braverman is a corresponding author. 
%Tel 1-403-2203956. Fax 1-403-282-5150.
}

\author[Braverman]{Elena Braverman}\ead{maelena@math.ucalgary.ca}, 
\author[Diblik]{Josef Dibl\'{i}k}\ead{diblik@feec.vutbr.cz},
\author[Rodkina,Diblik]{Alexandra Rodkina}\ead{alechkajm@yahoo.com},
\author[Diblik]{Zden\v{e}k \v{S}marda}\ead{smarda@feec.vutbr.cz}

\address[Braverman]{
Dept. of Math. \& Stats,
University of Calgary,
%2500 University Drive N.W., \
Calgary, %AB T2N 1N4, 
Canada
}

\address[Diblik]{Dept. of Math., Faculty of Electrical Engineering and Communication, Brno University of Technology, Brno, Czech Republic}
%\address[Smarda]{Dept. of Math., Brno University of Technology, Brno, Czech Republic}

%\author{Elena Braverman$^a$, Josef Dibl\'{i}k$^b$, 
%Alexandra Rodkina$^{b,c}$, Zden\v{e}k \v{S}marda$^b$}
%\date{last update: Feb 8, EB}

\address[Rodkina]
{Dept. of Math., the University of the West Indies,
%Mona, 
Kingston, Jamaica}


\begin{abstract}

Difference equations, such as a Ricker map, for an increased value of the parameter,  experience instability of the positive equilibrium and transition to deterministic chaos. To achieve stabilization, various methods can be applied. Proportional Feedback control suggests a proportional reduction of the state variable at every $k$th step. First, if $k \neq 1$, a cycle is stabilized rather than an equilibrium. Second, the equation can incorporate an additive noise term, describing the variability of the environment, as well as multiplicative noise corresponding to possible deviations in the control intensity. The present paper deals with both issues, it justifies a possibility  of getting a stable blurred $k$-cycle. Presented examples include
the Ricker model, as well as equations with unbounded $f$, such as the bobwhite quail population models. Though the theoretical results justify stabilization for either multiplicative or additive noise only, numerical simulations illustrate that a blurred cycle can be stabilized when both multiplicative and additive noises are involved.
 
%increasing for small $x$
%and vanishing at zero, we apply a noise-involving proportional feedback control every $k$th step
%$$
%x_{n+1}= \left\{ \begin{array}{ll} \displaystyle f\left( (\nu + \ell_1\chi_{m+1})x_n \right) + \ell_2\chi_{m+1},
%& n=mk, \\
%f(x_n), & n \neq mk,  \end{array} \right.~x_0>0, ~m,n \in {\mathbb N}_0, ~k\in {\mathbb N},~\nu \in (0,1].
%$$
%The purpose of getting a stable blurred $k$-cycle is achieved and illustrated with examples.
%Some generalizations are considered.
\end{abstract}

%{\bf AMS subject classification:} 39A50, 37H10, 34F05; Secondary:  39A30, 93D15, 93C55.

%{\bf Keywords:} 
%\subjclass{Primary:  39A50, 37H10, 34F05; Secondary:  39A30, 93D15, 93C55.}
% \keywords{Stochastic difference equations, proportional feedback control, stable $m$-cycles,
%population dynamics models, impulsive control}
%\email{maelena@ucalgary.ca}
%\email{alexandra.rodkina@uwimona.edu.jm}

%\thanks{
%E. Braverman is a corresponding author. e-mail maelena@ucalgary.ca.
%The first author is supported by NSERC grant RGPIN-2015-05976.}

%\begin{document}

\end{frontmatter}

%{\it last update: June 30, AR}

%{\footnotesize
%$^a$ Dept. of Math. \& Stats., University of Calgary, 2500 University Drive N.W. Calgary, AB,  T2N 1N4, Canada
%
%$^b$ Dept. of Math., Brno University of Technology, Brno, Czech Republic
%
%
%$^c$ Department of Mathematics,
%the University of the West Indies,
%Mona, Kingston, Jamaica}

%\centerline{\scshape Elena Braverman}
%\medskip
%{\footnotesize
%% please put the address of the first author
% \centerline{Dept. of Math. and Stats., University of Calgary}
% \centerline{2500 University Drive N.W. Calgary, AB,  T2N 1N4, Canada}}
%
%\medskip
%
%\centerline{\scshape  Alexandra Rodkina}
%\medskip
%{\footnotesize
%% please put the address of the first author
% \centerline{Department of Mathematics,
% the University of the West Indies,}
% \centerline{Mona, Kingston, Jamaica}}
%
%%\medskip

%\bigskip

%\centerline{(Communicated by Xiaoying Han)}

%\subjclass{Primary:  39A50, 37H10, 34F05; Secondary:  39A30, 93D15, 93C55.}
% \keywords{Stochastic difference equations, proportional feedback control,
%population dynamics models, Beverton-Holt equation}

\section{Introduction}
\label{sec:intr}

A difference equation 
\begin{equation}
\label{1}
x_{n+1}=f(x_n),\quad x_0>0, \quad n\in {\mathbb N}_0=\{0,1,2,\dots\}~,
\end{equation}
for a variety of maps $f$, for example, logistic or Ricker, can exhibit unstable and even chaotic behavior.
For unstable \eqref{1}, several control methods were developed in the literature, e.g. 
\cite{Dattani,TPC,gm,FL2010,LizPotsche14,uy99}.
These methods include Proportional Feedback (PF) control in the deterministic \cite{gm} and stochastic \cite{DCDSB2017} versions, 
Prediction-based control \cite{BKR,FL2010,LizPotsche14,uy99} and Target Oriented control \cite{Dattani,TPC}.  
Some of these methods were used to stabilize cycles rather than an equilibrium in \cite{NODY,Liz_CAMWA,LizPotsche14}.
Stochastic versions of these control methods, applied to stabilize a blurred equilibrium, were considered in \cite{BKR,DCDSB2017}. In addition, there are control methods where stabilization is achieved by noise only, see the recent papers 
\cite{GDEA_2019,Medv} and references therein.
In the present paper, we concentrate on a stochastic version of PF control, applied to stabilize blurred cycles. 

We consider the control by the proportional feedback (PF) method. 
This method, first introduced in \cite{gm}, 
involves reduction  of the state variable at each $k$-th step, 
$k\in {\mathbb N}$, when $n$ is divisible by $k$ ($n \mid k$),  
proportional to the size of the state variable $x_n$
\begin{equation}
\label{2}
x_{n+1}=f(\nu x_n), ~n \mid k, ~ x_{n+1}=f(x_n), ~n \not\,\mid k,
%x_{n+1}=\left\{ \begin{array}{ll} f(\nu x_n), & n \mid k , \\
%f(x_n), & n \not\,\mid k,  \end{array} \right. ~x_0>0,n\in {\mathbb N}_0,\nu \in (0,1].
\end{equation}
where $x_0>0$, $n\in {\mathbb N}_0$, %$m,n\in {\mathbb N}_0$, 
$\nu \in (0,1]$, $k\in {\mathbb N}$.

However, the reduction coefficient may involve a stochastic component, describing uncertainties in the control process, resulting in a  multiplicative noise
\begin{equation}
\label{3}
x_{n+1}= \left\{ \begin{array}{ll} \displaystyle f\left( (\nu + \ell_1\chi_{n+1})x_n \right),
& n \mid k, \\
f(x_n), & n \not\,\mid k, \end{array} \right.
%~m,n \in {\mathbb N}_0,~ x_0>0, ~ \nu \in (0,1],~ k\in {\mathbb N}.
\end{equation}
$x_0>0$, $n\in {\mathbb N}_0$, 
$\nu \in (0,1]$, $k\in {\mathbb N}$.
We can also consider the case when the reduction coefficient is deterministic but there are random fluctuations of $x_n$
at the control step, describing variability of the environment
\begin{equation}
\label{4}
x_{n+1}= \left\{ \begin{array}{ll} \max \left\{ f(\nu x_n) +\ell_2\chi_{n+1}, 0 \right\},
& n \mid k, \\
f(x_n),%+\ell\chi_{n+1}, 
& n \not\,\mid k,  \end{array} \right.
\end{equation}
$x_0>0$, $n \in {\mathbb N}_0$, $\nu \in (0,1]$, $k\in {\mathbb N}$.
Here $\chi_{n+1} \in [-1,1]$ is the bounded random variable, while $\ell_j$, $j=1,2$ describes the bound of the noise. 
%While \eqref{3} accounts for possible fluctuations of harvesting effort, model \eqref{4} considers a random deduction at %each control step, 
%which can describe either pollution and disturbance associated with a harvesting event, 
%as well as consumption or sampling applied by controllers and independent of the current population size.
%On a positive side, harvesters may cause increase by unintentionally fostering immigration or boosting available food %supply. 

The deterministic version of cycle stabilization by PF control was justified in \cite{NODY}.
Stabilization of a positive equilibrium with PF method shifts an equilibrium closer to zero and
is achieved in an interval $\nu \in (\alpha,\beta) \subset (0,1)$. 
For smaller values of $\nu$, zero becomes the only 
stable equilibrium, for higher values, a positive equilibrium still can be unstable.
When we applied PF control on each $k$th step \cite{NODY}, it led to %construction of 
an asymptotically stable $k$-cycle, 
with all the values between zero and a positive equilibrium. Here we construct a stochastic analogue of this process.
%%%%%%%%%%%%%%%%%%%%%%
Stabilization of stochastic equations with proportional feedback was recently explored in the continuous case 
\cite{IoM}, as well as the idea of periodic controls \cite{YSZL}.
%%%%%%%%%%%%%%%%%%%%%%%

The paper is organized as follows.
In Section~\ref{sec:prelim}, we introduce main assumptions %for function $f$ and stochastic perturbations, 
and discuss properties of a $k$-iteration of function $f$. Section \ref{sec:multi} contains results on the existence of a blurred $k$-cycle in the presence of stochastic  multiplicative perturbations of the 
control parameter $\nu$ when the level of noise $\ell$ is small, while Section \ref{sec:add} deals with 
the controlled equation for additive stochastic perturbations. Section \ref{sec:ex} contains examples 
with computer simulations illustrating the results of the paper, along with some generalizations. In particular, 
a modification of PF method ``centered'' at an unstable equilibrium $K$ instead of zero, is developed and applied 
to construct a blurred $k$-cycle in the neighborhood of $K$, when both stochastic,  multiplicative and 
additive perturbations, are present.

%%%%%%%%%%%%%%%%%%%%%%%%
\section{Definitions  and  Auxiliary Results}
\label{sec:prelim}

In this paper, we impose an  assumption on the map $f$ %which describes its behaviour 
in a right neighbourhood  of zero.

\begin{assumption}
\label{as:slope}
\emph{The function $f : [0, \infty) \rightarrow [0, \infty)$ is continuous, $f(0)=0$, and
there is a real number $b>0$ such
that $f(x)$ is strictly monotone increasing, while the function $f(x)/x$ is strictly monotone
decreasing on $(0,b]$, $f(b) > b$, while $f(b)/b> f(x)/x$ for any $x \in (b, \infty)$.}
\end{assumption}

\begin{remark}
\label{rem:bb1}
 Note that, once  Assumption \ref{as:slope} holds for a certain $b>0$, it is also satisfied for any $b_1\in (0, b]$.
\end{remark}

Many functions in \eqref{1} used in applications satisfy Assumption~\ref{as:slope}, see \cite{Thieme} and examples below. We truncate values 
at zero, when necessary, to satisfy $f : [0, \infty) \rightarrow [0, \infty)$, which is a common practice \cite{schreiber}.
Examples include the Ricker model
\begin{equation}
\label{eq:ricker}
x_{n+1} = f_1(x_n) = x_n e^{r(1-x_n)}
\end{equation}
for $r>1$,  with any $b \leq 1/r$,
the logistic model (truncated at zero)
%\begin{equation*}
%\label{eq:logistic}
$\displaystyle x_{n+1} = f_2(x_n) = \max\left\{ rx_n (1-x_n), 0 \right\}$
%\end{equation*}
for $r>2$, with $b \leq 1/2$.  
%%% omit for Automatica
%The Maynard Smith model \cite{Thieme}
%$\displaystyle
%x_{n+1} = f_3(x_n) = \frac{Ax_n}{1+Bx_n^{\gamma}}$, where $A,B>0$, $\gamma>1$,
%and another modification of the Beverton-Holt equation
%$\displaystyle x_{n+1} = f_4(x_n) = \frac{Ax_n}{(1+Bx_n)^{\gamma}}$, $A,B>0$, $\gamma>1$
%satisfy Assumption~\ref{as:slope}, with certain parameters, 
%leading to the existence of a positive equilibrium exceeding the  unique maximum point $x_{\max}$, with 
%$b<x_{\max}$. 
%%% 
In these maps, $f_i$ are unimodal, increasing on $[0,x_{\max}]$ and decreasing on $[x_{\max},\infty)$, with the only critical point on 
$[0,\infty)$, which is a global maximum.
However, Assumption~\ref{as:slope} can hold for functions which have more than one critical point, for example,
for the map developed in \cite{Milton} to describe the growth of the bobwhite quail population
\begin{equation}
\label{eq:milton} 
f_3(x) = x \left( A + \frac{B}{1+x^{\gamma}} \right), \quad A,B > 0, \quad \gamma>1,
\end{equation}
which, generally, has two critical points, first a local maximum, then a global minimum, then increases, and $f_3(x) \to \infty$ as $x \to \infty$.

We denote by  $(\Omega, {\mathcal{F}},  (\mathcal{F}_m)_{m \in \mathbb{N}}, {\mathbb{P}})$   a complete
filtered probability space, $\chi:=(\chi_m)_{m\in\mathbb{N}}$ is  a sequence of independent random variables with the zero mean. 
The filtration $(\mathcal{F}_m)_{m \in \mathbb{N}}$ is naturally generated by  the sequence $(\chi_m)_{m\in\mathbb{N}}$, i.e.
$\mathcal{F}_{m} = \sigma \left\{\chi_{1},  \dots, \chi_{m}\right\}$.
The standard abbreviation ``a.s." is used for both ``almost sure" or ``almost surely" with respect
to the fixed probability measure $\mathbb P$  throughout the text.
A detailed discussion of stochastic concepts and notation can be found in \cite{Shiryaev96}.
%In the present paper, 
We consider \eqref{3} and \eqref{4},
where the sequence $(\chi_m)_{m \in {\mathbb N}}$ satisfies the following condition.

\begin{assumption}
\label{as:chibound}
\emph{$(\chi_m)_{m \in {\mathbb N}}$ is a sequence of independent and identically  distributed
continuous random variables, with the density function $\phi(x)$ such that
$\phi(x)>0$ for $x\in [-1, 1]$ and $\phi(x) = 0$ for $x\notin [-1, 1]$.}
\end{assumption}

\begin{remark}
\label{rem:chibound}
In fact, Assumption~\ref{as:chibound} can be relaxed to the condition
$\displaystyle {\mathbb{P}} \left\{ \chi \in [1-\varepsilon,1] \right\} >0$ 
for any $\varepsilon>0$, which would allow to include discrete 
distributions, where $\displaystyle {\mathbb{P}} \left\{ \chi =1 \right\} >0$.
% with $\chi$ taking a value of one with a positive probability.
\end{remark}

In numerical simulations, we also consider the combination of \eqref{3} and \eqref{4}
\begin{equation}
\label{5_new}
\vspace{-0.3cm} x_{n+1}= \left\{ \begin{array}{ll} \displaystyle f\left( (\nu + \ell_1\chi_{n+1})x_n \right) + \ell_2\chi_{n+1},
& n \mid k, \\
f(x_n), & n \not\,\mid k,  \end{array} \right.
\end{equation}
$x_0>0$, $n \in {\mathbb N}_0$, $k\in {\mathbb N}$, $\nu \in (0,1]$.
%$$x_0>0, \quad n\in \mathbb N, \quad \nu \in (0,1], \quad k\in {\mathbb N}.$$

Let us start with some auxiliary results on $f^k(x)=f(f^{k-1}(x))$ and $g(x):=f^k(\nu x)$ for any $\nu \in (0,1]$,
where Assumption~\ref{as:slope} holds.
%the function $f(\nu x)$, for any $\nu \in (0,1]$, where Assumption~\ref{as:slope} holds. 
Obviously $f:[0, b] \to [0, f(b)]$ is  increasing and continuous,
and there is an increasing  and continuous inverse  function
$f^{-1}:[0, f(b)]\to [0, b]$. 
As $f(b) > b$ and $f(x)/x$ is decreasing on $(0,b]$ by Assumption \ref{as:slope}, $f(x)>x$ for $x\in [0,b]$, and also $f$ is increasing. 
Thus $f^{-1}(b):[0,f(b)]\to [0,b]$ is well defined, and $f^{-1}(b) \in (0,b)$. 
Evidently $f^2:[0,f^{-1}(b)] \to [0,f(b)]$ is continuous and increasing, since $f$ is increasing on $[0,b]$, 
and $f^2(x) \in [0,f(b)]$ for $x\in [0,f^{-1}(b)]$. 
Therefore $f^{-2}:[0,f(b)] \to [0,f^{-1}(b)]$ is also well defined and increasing. 
Similarly, $f^{-k}:[0,f(b)] \to [0,f^{1-k}(b)]$ exists and is increasing for any $k \in {\mathbb N}$.  
Denote
\begin{equation}
\label{def_bk}
b_j := f^{1-j} (b), \quad j \in {\mathbb N},\quad j \neq 1, \quad b_1=b,
\end{equation}
then $f(b_{j+1})=b_j$, $j=1,2, \dots, k$, and
\begin{equation}
\label{bk_ineq}
b=b_1 > b_2 > \dots > b_k  >0.
\end{equation}

\begin{lemma}
\label{iterate}  
If $f$ satisfies Assumption \ref{as:slope}, this assumption also holds for $f^k$ with
$b_k$ instead of $b$, where $b_k$ is defined in \eqref{def_bk}. 
\end{lemma}
%\begin{proof}
{\bf Proof.}
The function $f:[0,b] \to [0,f(b)]$ is continuous and monotone increasing, so is $f^k:[0,b_k]=[0,f^{1-k}(b)] \to [0,f(b)]$. Next, let us prove that 
$f^k(x)/x$ is monotone decreasing on $[0,b_k]$. Let $0<x_1<x_2 \leq b_k$. 
Then $f(x_1) \leq b_{k-1}$, \dots, $f^j(x_1)\leq b_{k-j}$, $j=1, \dots, k-1$. Since $f(x)/x$ is decreasing on $[0,b]$, while $f$ is increasing,
$f(x_1)/x_1 > f(x_2)/x_2$, $f(f(x_1))/f(x_1)>f(f(x_2))/f(x_2)$, \dots, $f^k(x_1)/f^{k-1}(x_1)>f^k(x_2)/f^{k-1}(x_2)$
and
$$
\frac{f^k(x_1)}{x_1}= \frac{f^k(x_1)}{f^{k-1}(x_1)} \dots \frac{f^2(x_1)}{f(x_1)} \frac{f(x_1)}{x_1}
$$
$$
> \frac{f^k(x_2)}{f^{k-1}(x_2)} \dots \frac{f^2(x_2)}{f(x_2)} \frac{f(x_2)}{x_2}
=\frac{f^k(x_2)}{x_2}.
$$
Also, $f(0)=0$ implies $f^k(0)=0$, and $f(b) > b$, by \eqref{bk_ineq},  yields that 
$\displaystyle f^k(b_k)=f^k\left( f^{1-k} (b)   \right) = f(b) > b > b_k.$

Finally, let us justify that $f^k(x)/x<f^k(b_k)/b_k$ for any $x>b_k$ by induction. 
For $k=1$,  $f(x)/x<f(b)/b$ follows from Assumption~\ref{as:slope}. 

For $k=2$ and $x>f^{-1}(b)=b_2$, we consider two possible cases: $f(x)<b$ and $f(x) \geq b$.
In the former case, $f(f(x)) < f(b)$, as $f$ increases on $[0,b]$, and 
$$
\frac{f^2(x)}{x} < \frac{f(b)}{x} < \frac{f(b)}{b_2}=\frac{f^2(b_2)}{b_2}.
$$
For $f(x) \geq b$, by Assumption~\ref{as:slope}, $f(x)/x<f(b_2)/b_2$ for any $x>b_2$, as $b_2 < b$,
and $f(f(x))/f(x) \leq f(b)/b$.
Thus
$$\frac{f(f(x))}{x} = \frac{f(f(x))}{f(x)} \frac{f(x)}{x} < \frac{f(b)}{b} \frac{f(b_2)}{b_2}$$ $$=
\frac{f(b)}{b}\frac{b}{b_2}=\frac{f(b)}{b_2}=\frac{f^2(b_2)}{b_2}.
$$
Next, let us proceed to the induction step. 
Assume $\displaystyle \frac{f^n(x)}{x}<\frac{f^n(b_n)}{b_n}= \frac{f(b)}{b_n}$ for any $x>b_n$. Consider $x>b_{n+1}$.
Then either $f^n(x) < b$ or $f^n(x) \geq  b$. In the former case $f^n(x) < b$, we have $f(f^n(x)) < f(b)$ due to monotonicity of $f$ on 
$[0,b]$ and
$$
\frac{f^{n+1}(x)}{x} < \frac{f(b)}{x} < \frac{f(b)}{b_{n+1}}=\frac{f^{n+1}(b_{n+1})}{b_{n+1}}.
$$
In the latter case $f^n(x) \geq  b$ we get
$$
\frac{f^{n+1}(x)}{x} = 
\frac{f(f^{n}(x))}{f^n(x)} \frac{f^n(x)}{x} < 
\frac{f(b)}{b} \frac{f^n(b_{n+1})}{b_{n+1}}$$ $$= \frac{f^{n+1}(b_{n+1})}{b} \frac{b}{b_{n+1}} =  
\frac{f^{n+1}(b_{n+1})}{b_{n+1}},
$$
where in the inequality we used $\displaystyle \frac{f^n(x)}{x}<\frac{f^n(b_{n+1})}{b_{n+1}}$ for any $x>b_{n+1}$
by the induction assumption. Also, $f(u)/u \leq f(b)/b$ for any $u=f^n(x) \geq b$ by Assumption~\ref{as:slope}, while equalities applied notation 
\eqref{def_bk}.
%This concludes the proof.
%\bigskip
%\end{proof}
\qed

Define the function $\Psi_k$ as
\begin{equation}
\label{def:Psi}
\Psi_k(x):=\frac {x}{f^k(x)}, \quad x\in (0, b_k), \quad k \in {\mathbb N},
\end{equation}
and formally introduce the limit
\begin{equation}
\label{def:psi0}
\Psi_k(0):= \lim\limits_{x \to 0^+} \frac {x}{f^k(x)}.
\end{equation}

\begin{lemma}
\label{lem:psi-1}
Let Assumption \ref{as:slope} hold, $k \in {\mathbb N}$ 
and $\Psi_k$ be defined as in \eqref{def:Psi}, \eqref{def:psi0}.  
Then 
\begin{enumerate}
\item $\Psi_k: (0, b_k)\to \biggl( \Psi_k(0),  \Psi_k(b_k)\biggr)$, \\ $\Psi_k^{-1}: \biggl(\Psi_k(0),  
\Psi_k(b_k)\biggr) 
\to (0, b_k)$;
\item $0\le \Psi_k(0)<\Psi_k(b_k)<1$;
\item  both $\Psi_k$ and its inverse $\Psi_k^{-1}$ are increasing and continuous on their domains.
\end{enumerate}
\end{lemma}
{\bf Proof.}
By Lemma~\ref{iterate}, the function $\Psi_k$ defined in \eqref{def:Psi} 
%\begin{equation}
%\label{def:Psi}
%\Psi_k(z):=\frac {z}{f^k(z)}, \quad z\in (0, b_k), \quad k \in {\mathbb N}
%\end{equation}
is increasing, continuous and hence has a unique
inverse function on $(0, b_k)$.
Following Assumption \ref{as:slope}, we notice that the limit $\displaystyle \lim_{x \to 0^+} \frac{f(x)}{x}$
%\begin{equation*}
%\label{def:lim}
% \lim_{x \to 0^+} \frac{f(x)}{x}
%\end{equation*}
exists (finite or infinite), is positive and greater than 1, since $f(x)/x$ is
decreasing on $(0,b_k)$ and $f^k(b_k) > b_k$.  Note that
$\displaystyle
\frac{f^k(x)}{x}=\frac 1{\Psi_k(x)}$
and 
%$$
$\displaystyle \lim_{x \to 0^+} \frac{f^k(x)}{x}=\lim_{x \to 0^+} \frac 1{\Psi_k(x)}$,
%$$
where $\displaystyle \lim_{x \to 0^+} \frac{1}{\Psi_k(x)}=0$ if
$\displaystyle \lim_{x \to 0^+} \frac{f^k(x)}{x}=+\infty$. Thus \eqref{def:psi0} is well defined,
and Part (1) is valid. 
%As mentioned above
Also, $\Psi_k$ and its inverse are continuous monotone increasing in their domains and by Lemma~\ref{iterate}, $\Psi_k(b_k)<1$, which 
implies Parts (2) and (3).
\qed
%\end{proof}
%\bigskip 

To apply known results from \cite{DCDSB2017}, for each point $x^*\in (0, f(b))$,  we are looking for the control parameter $\nu=\nu(x^*)\in (0, 1)$  such that 
$x^*$ is a fixed  point of the function $g(x):=f^k(\nu x)$.  We  recall from \eqref{def_bk} that $b_1=b$ and introduce
\begin{equation}
\label{def:nux*}
\begin{array}{l}
\hat{x}= f^{-k}(x^*),~\nu=\nu(x^*):=\Psi_k(f^{-k}(x^*))\\ \nu(x^*)=\Psi_k(\hat{x}),~ 
\hat{x}=\nu(x^*)x^*=\Psi_k^{-1}(\nu). \end{array}
\end{equation}

\begin{lemma}
\label{lem:fixg}
Let Assumption \ref{as:slope} hold, $k \in {\mathbb N}$ and $x^*\in (0, f(b))$.
The function $\nu(x^*)$ defined in \eqref{def:nux*} satisfies the following conditions:
\begin{enumerate}
%\end{enumerate}
%\begin{enumerate}
%(i) 
\item
$x^*$ is a fixed point of $g(x):=f^k(\nu x)$, i.e. $f^k\bigl(\nu(x^*)x^*\bigr)=x^*$, $\nu(x^*)=\Psi_k(\hat{x})$,
$\hat{x} \in (0,b_k)$; 
%(ii) 
\item
$\nu(x^*)\in \left(\Psi_k(0), \,  \Psi_k(b_k)\right)\subset (0, 1)$; %\vspace{1.5mm}
%(iii) 
\item
$\nu(x^*)$ is an increasing function of $x^*$ on $(0, f(b) )$.
\end{enumerate}
\end{lemma}
%\begin{proof}
{\bf Proof.} 
(1). Let $x^*\in (0,f(b))$, then $f^{-k}(x^*) \in (0, b_{k})$, thus $\Psi_k (f^{-k}(x^*))$ is well defined  and
$
\nu(x^*)=\Psi_k (f^{-k}(x^*))$ \\ $\displaystyle
=\frac{f^{-k}(x^*)}{f^k(f^{-k}(x^*))}=\frac{f^{-k}(x^*)}{x^*}=\frac{\hat{x}}{f^k(\hat{x})}
=\Psi_k(\hat{x}),
$
hence $\hat{x} =\nu(x^*)x^* \in (0,b_k)$ and
$$
f^k\bigl(\nu(x^*)x^*\bigr)=f^k\left(\frac{f^{-k}(x^*)}{x^*} x^*\right)=f^k\left(f^{-k}(x^*)\right)=x^*.
$$
(2). We have $x^*\in (0, f(b))$ and
$\hat{x}=f^{-k}(x^*)\in (0, b_{k})$. Thus Lemma~\ref{lem:psi-1}, Part 1 implies
$
\nu(x^*)\in \left(\Psi_k(0), \,  \Psi_k(b_k)\right)\subset (0, 1).
$
\\
(3). By Lemma~\ref{lem:psi-1} and Assumption~\ref{as:slope},
for any $k \in {\mathbb N}$, both $\Psi_k$ and $f^{-k}$ are increasing functions on $(0, b_k)$ and $(0, f(b))$, 
respectively. 
Therefore $\nu(x^*)=\Psi_k(f^{-k}(x^*))$  is increasing as a function of $x^*$ on $(0, f(b))$.%, which concludes the proof.
\qed
%\medskip

%\end{proof}

\section{Multiplicative perturbations}
\label{sec:multi}

Consider the deterministic PF with variable intensity $\nu_m \in (0,1]$, 
applied at each $k$-th step, for a fixed $k \in \mathbb N$,
\begin{equation}
\label{eq:var_imp}
%\vspace{-0.3cm}
x_{n+1}= \left\{  \begin{array}{ll} f(\nu_n x_n), & n \mid k, \\
f(x_n), & n \not\,\mid k,  \end{array} \right.  x_0>0, n\in {\mathbb N}_0.
\end{equation}
%$m,n\in {\mathbb N}_0$. 
Investigation of \eqref{eq:var_imp} will allow to analyze corresponding stochastic equation \eqref{3} with a multiplicative noise.  
For each $x^*\in \bigl( 0, f(b)\bigr)$,  we establish the control
$\nu=\nu(x^*)$ and define an interval such that
a solution of \eqref{3} remains in this interval, once the level of noise $\ell$ is small enough.
This method was applied, for instance, in \cite{BKR}. 

Further, we apply the result obtained in \cite{DCDSB2017} for 
\begin{equation}
\label{eq:var}
z_{m+1}= g\left( \nu_m z_m \right) = f^{k}\left( \nu_m z_m \right), ~ z_0>0,~ m\in \mathbb N,
\end{equation}
to explore stochastic equation \eqref{3} with a multiplicative noise.

For any $\mu_1$, $\mu_2$ such that
\begin{equation}  
\label{eq:mu_bounds}
\Psi_k(0) < \mu_1<   \mu_2 < \Psi_k(b_k), 
\end{equation}
we define
\begin{equation}   
\label{eq:defys}
y_1:= \Psi_k^{-1}(\mu_1), \quad y_2 :=\Psi^{-1}_k(\mu_2).
\end{equation}

\begin{lemma} \cite[Lemma 3.1]{DCDSB2017}
\label{lem:multi}
Let Assumption~\ref{as:slope}  hold for $f^k$, $k\in {\mathbb N}$,  $\mu_1$ and $\mu_2$ satisfy \eqref{eq:mu_bounds} 
and, for each $m \in {\mathbb N}$,
\begin{equation}
\label{eq:nu_bounds}
\nu_m \in [\mu_1,\mu_2].
\end{equation}
Then, for any $z_0>0$ and $\varepsilon$, $\displaystyle 0<\varepsilon < \min\left\{ y_1, b_k-y_2\right\}$, where
$y_1$,$y_2$ are defined in \eqref{eq:defys},
%$$
%0<\varepsilon < \min\left\{  \Psi_k^{-1}(\mu_1), b_k-\Psi_k^{-1}(\mu_2)  \right\},
%$$
there is $m_0=m_0(x_0,\varepsilon)$, $m_0 \in {\mathbb N}$, such that
the solution $z_n$ of equation \eqref{eq:var} for any $m \geq m_0$ satisfies
\begin{equation}
\label{eq:z_bounds}
\nu_m z_m \in \left(y_1-\varepsilon, y_2+\varepsilon \right).
%\nu_m z_m \in \left(\Psi_k^{-1}(\mu_1)-\varepsilon, \Psi_k^{-1}(\mu_2)+\varepsilon\right).
\end{equation}
\end{lemma}

\begin{remark}
Lemma~\ref{lem:multi} actually states (see its proof in \cite{DCDSB2017}) that, 
for a prescribed $k \in {\mathbb N}$, for a small $\varepsilon>0$,
once a solution of \eqref{eq:var_imp} satisfies
$\nu_{km} x_{km}\in (y_1-\varepsilon, y_2+\varepsilon)$, $m \in {\mathbb N}$,
all the subsequent $k$-iterates $\nu_{m+j} x_{(m+j)k}$, $j \in {\mathbb N}$, are also in this interval.
This is also true for the results based on Lemma~\ref{lem:multi}, in particular,
for Lemma~\ref{lem:multistoch} and Theorem~\ref{thm:mainmult}.
\end{remark}

\begin{lemma} 
\label{lem:multi_cycle}
Let Assumption~\ref{as:slope}  hold,   $\mu_1$, $\mu_2$ satisfy \eqref{eq:mu_bounds} 
and, for each $m \in {\mathbb N}$, \eqref{eq:nu_bounds} be fulfilled.
For any $x_0>0$ and $\varepsilon>0$, 
%$\displaystyle \varepsilon < \Psi_k^{-1}(\mu_1)$, $\displaystyle \varepsilon < b_k-\Phi_k^{-1}(\mu_2)$.
there is $m_0=m_0 (x_0, \varepsilon) \in {\mathbb N}$ such that for $m\geq m_0$,
the solution of \eqref{eq:var_imp} satisfies
\begin{equation}
\label{eq:sol_mul_bounds}
x_{mk+j}\in \left(f^j(y_1)-\varepsilon, f^j(y_2)+\varepsilon \right), \quad j=1, \dots, k,
\end{equation}
where $y_1$ and $y_2$ are defined in \eqref{eq:defys}.
\end{lemma}
{\bf Proof.}
%\begin{proof} 
Note that $y_1,y_2\in (0,b_k)$ and $f^j$ are  continuous and monotone increasing on $(0,b_k)$ for $j=0, \dots, k-1$.
Therefore for any $\varepsilon>0$ there is an $\varepsilon_1>0$ such that for $j=1, \dots, k$,
%\begin{equation}
%\label{eq:epsilon}  
%\vspace{-0.2cm}
%\begin{array}{l}
$\displaystyle u \in \left(y_1-\varepsilon_1, y_2+\varepsilon_1 \right) ~\Rightarrow $ \\ %\\ 
$\displaystyle f^j(u) \in \left(f^j(y_1)-\varepsilon, 
f^j(y_2) + \varepsilon \right)$. %, \end{array} ~ j=1, \dots, k.
%\end{equation}
Choose $z_0=x_0$ and
$\displaystyle
\varepsilon_2 < \min \left\{  y_1,b_k-y_2,\varepsilon_1 \right\}$
instead of $\varepsilon$  in Lemma~\ref{lem:multi}. Then for $m>m_0$,  by \eqref{eq:z_bounds}, $\nu_{mk}x_{mk} \in (y_1-\varepsilon_2,
y_2+\varepsilon_2)$. Since, by the above implication, %\eqref{eq:epsilon}, 
\\
$
x_{mk+1}=f\left(\nu_{mk}x_{mk} \right) \in \left(f(y_1)-\varepsilon,f(y_2) + \varepsilon \right),
\dots$, \\ $ x_{mk+k}=f^k\left(\nu_{mk}x_{mk} \right) \in \left(f^k(y_1)-\varepsilon,f^k(y_2) + \varepsilon \right),
$ \\
this implies \eqref{eq:sol_mul_bounds}.% and concludes the proof.
\qed
%\end{proof}
%\medskip

Let us proceed to stochastic equation~\eqref{3}.

We start with an auxiliary result which follows from Lemma~\ref{lem:multi_cycle}.

\begin{lemma}
\label{lem:multistoch}
Let $k \in {\mathbb N}$ be fixed, Assumptions~\ref{as:slope} and~\ref{as:chibound} hold,
$\Psi_k$ be defined in \eqref{def:Psi},  $x^*\in (0, f(b))$, $\nu=\nu(x^*)$ be as in \eqref{def:nux*},
and %$\ell\in \mathbb R$ satisfy the inequality
\begin{equation}
\label{def:alphal}
\ell  \in \bigl( 0, \min\bigl\{\Psi_k(b_k)-\nu, \, \, \nu -  \Psi_k(0) \bigr\} \bigr), 
\end{equation}
%and
\begin{equation}
\label{y_bounds}
\underline{y} :=  \Psi_k^{-1}(\nu-\ell),~ \overline{y} :=   \Psi_k^{-1}(\nu+\ell),~ 
0< \underline{y} < \overline{y} < b_k. 
\end{equation}
Let $x_n$ be a solution to equation~\eqref{3} with $\nu$,$\ell$ satisfying \eqref{def:nux*} and
\eqref{def:alphal}, respectively. 

Then, for any $\varepsilon>0$ there is a $m_0=m_0(\varepsilon, x^*, x_0) \in {\mathbb N}$ such that, 
for all $m\ge m_0$, $m \in {\mathbb N}$,
$$
x_{mk+j}\in \left(f^j(\underline{y})-\varepsilon, f^j(\overline{y})+\varepsilon \right), \quad j=1, \dots, k, \mbox{~~a.s.}
$$
\end{lemma}
{\bf Proof.}
Since $x^*< f(b_k)<f(b)$, Lemma \ref{lem:fixg} implies
$\displaystyle \nu(x^*)=\Psi_k (f^{-k}(x^*))\in \left(\Psi_k(0), \, \Psi_k(b_k)\right)$.
Thus the right segment bound $\nu - \Psi_k(0)$  in  \eqref{def:alphal} is positive.  By 
Assumption \ref {as:chibound} we have, a.s.,
$$\nu_{m}=\nu+\ell \chi_{m+1} \leq \nu+\ell, \quad \nu_m=\nu+\ell \chi_{m+1} 
\geq \nu - \ell
$$
and $\displaystyle \nu_m=\nu+\ell \chi_{mk+1} \geq \nu - \ell$, thus $\nu_m \in [\nu-\ell,\nu+\ell]$, a.s.
Let 
 $\mu_1:= \nu-\ell$, $\mu_2:= \nu+\ell$. 
With $\nu$ as in \eqref{def:nux*} and
$\ell$ satisfying  \eqref{def:alphal}, we have
\[
\Psi_k(0)-\ell<\mu_1<\mu_2<\Psi_k(b_k)+\ell, 
\]
%then, by 
then Lemma~\ref{lem:multi_cycle} implies the statement of the lemma.
\qed
%which concludes the proof. 
%\end{proof}
%\medskip

Lemma \ref{lem:multistoch} leads to the main result of this section, which states that for 
any $k \in {\mathbb N}$ and
$x^*\in (0, f(b))$, we can find a control $\nu$ and a noise level $\ell$, such that  the solution
eventually reaches some neighbourhood of a $k$-cycle, a.s., and stays there.

\begin{theorem}
\label{thm:mainmult}
Let Assumptions~\ref{as:slope} and \ref{as:chibound} hold, $\Psi_k$  be defined in 
\eqref{def:Psi}, \eqref{def:psi0},
$x^*\in (0, f(b))$ be an arbitrary point,
$\nu=\nu(x^*)$ be denoted in \eqref{def:nux*}, 
$\underline{y}$ and $\overline{y}$ be defined in \eqref{y_bounds},
$x_0>0$ and $\ell\in \mathbb R$ satisfy inequality
\eqref{def:alphal}. Then for the solution $x_n$ of equation~\eqref{3}, the following statements hold. 
\\
%\begin{enumerate}
(i) For each $\varepsilon>0$ there exists a nonrandom $m_0=m_0(\varepsilon, x^*, x_0)\in \mathbb N$ 
such that, for all $m\geq m_0$,
\\
$ \displaystyle
x_{mk+j} \in \left( f^j(\underline{y})-\varepsilon, f^j(\overline{y})+\varepsilon \right), \quad j=1, \dots, k, \mbox{~~a.s.}
$
%where $\underline{y}$ and $\overline{y}$ are defined in \eqref{y_bounds}.
%\begin{equation*}
%\label{eq:intmultn_cycle}
%x_n\in\left (\frac{\Phi^{-1}(\nu-\ell)} {\nu+\ell}-\delta, \,\,  \frac{\Phi^{-1}(\nu+\ell)} {\nu-\ell}+\delta\right).
%\end{equation*}
%\\
(ii)
$\displaystyle
\liminf_{m \to \infty} x_{mk+j} \geq  f^j(\underline{y}),
~ \limsup_{m \to \infty} x_{mk+j} \leq f^j(\overline{y}), ~j=1, \dots, k, 
~ \text{a.s.}
$
%\end{enumerate}
\end{theorem}
%\medskip
{\bf Proof.}
%\begin{proof}
Note that from condition \eqref{def:alphal} we have $\nu>\ell$.
%Fix $\delta>0$ and take some $\varepsilon>0$ satisfying
%\[
%\varepsilon<\delta (\nu-\ell).
%\]
By Lemma~\ref{lem:multistoch}, for any $x_0>0$ and $\varepsilon>0$, there is $m_0=m_0(\varepsilon)\in \mathbb N$ such that, a.s.,
$
x_{mk+j} > f^j(\underline{y})-\varepsilon$, $x_{mk+j}< f^j(\overline{y})+\varepsilon$, $m \geq m_0$, $j=1, \dots, k$,
which immediately implies (i).

Choosing a sequence of $\varepsilon_m=\frac{1}{m}$, $m \in \mathbb N$ in (i), 
we deduce (ii).
\qed
%\end{proof}
%\bigskip

Next, let us assume that the level of noise can be chosen arbitrarily small.
Theorem~\ref{prop:liml0} below 
%deals with the situation when the noise level $\ell$ can be chosen arbitrarily small.  It 
confirms the intuitive feeling that, as the noise level $\ell$ is getting smaller, the solution of  stochastic equation \eqref{3}  
behaves similarly to the solution of corresponding deterministic equation \eqref{2}  in terms of  approaching its stable cycle $\{ 
f^j(\hat{x})\}$, $j=1, \dots, k$, where $\hat{x}$ is defined in \eqref{def:nux*}.

\begin{theorem}
\label{prop:liml0}
Let Assumptions~\ref{as:slope} and~\ref{as:chibound} hold, $k\in {\mathbb N}$ be fixed, $\hat{x} \in (0, b_k)$ be an arbitrary point, 
$x^*=f^k(\hat{x})$, $\nu=\nu(x^*)$ be defined as in \eqref{def:nux*}, and $x_0>0$.
Then, for any   $\varepsilon>0$,  there exists the level of noise $\ell(\varepsilon)>0$ such that
for each $\ell <\ell (\varepsilon)$, there is a nonrandom  $m_1=m_1(\varepsilon, \ell, \hat x, x_0)$ such that the solution
$x$  of equation~\eqref{3} satisfies $x_{mk+j} \in (f^j(\hat{x})-\varepsilon, f^j(\hat{x})+\varepsilon)$, $j=1, \dots, k$ 
for $m\ge m_1$, a.s.
\end{theorem}
{\bf Proof.} First of all, from monotonicity of $f^k$ notice that the map $x^*=f^k(\hat{x})$ is one-to-one, 
and an arbitrary $x^* \in (0,f(b))$ corresponds to a certain $\hat{x} \in (0, b_k)$.
Next, 
by continuity of all $f^j$, %$j=1, \dots, k$, 
for any $\nu=\nu(x^*)$ defined as in \eqref{def:nux*},
there is a $\delta>0$ such that
\begin{equation}  
\label{add_th_3_5}
|x-\hat{x}| < \delta  \Rightarrow  \left| f^j(x)-f^j(x^*) \right| < \frac{\varepsilon}{2},~ j=1, \dots, k.
\end{equation}
Also, from the choice of $\nu$ in \eqref{def:nux*} and continuity of $\Psi_k$, % on a chosen interval,
there is $\ell(\varepsilon)>0$ such that for $\ell <\ell (\delta)$, %we have
$$
\left| \underline{y} - \hat{x} \right|< \delta, \quad  \left| \overline{y} - \hat{x} \right| < \delta,
$$
since $\underline{y}$ and $\overline{y}$ defined in \eqref{y_bounds} continuously depend on $\ell$.
Thus, by \eqref{add_th_3_5},
\begin{equation}
\label{add1_iter}
\left| f^j(\hat{x}) - f^j(\underline{y}) \right| < \frac{\varepsilon}{2},   \left| 
f^j(\overline{y})-f^j(\hat{x}) \right|
 < \frac{\varepsilon}{2}, 
\end{equation}
$j=1, \dots, k$.
Next, let us apply Theorem~\ref{thm:mainmult}, Part (i), with $\frac{\varepsilon}{2}$ instead of $\varepsilon$. 
Then, $\forall m\geq m_0$, $j=1, \dots, k$,  a.s.,
\begin{equation}
\label{add2}
x_{mk+j} \in \left( f^j(\underline{y})-\frac{\varepsilon}{2}, f^j(\overline{y})+\frac{\varepsilon}{2} \right).
\end{equation}
In view of \eqref{add1_iter} and \eqref{add2}, $\displaystyle x_{mk+j}> 
f^j(\hat{x})-\frac{\varepsilon}{2}-\frac{\varepsilon}{2}=f^j(\hat{x})-\varepsilon$
and $\displaystyle x_{mk+j}< f^j(\hat{x})+\frac{\varepsilon}{2}+\frac{\varepsilon}{2}=f^j(\hat{x})+\varepsilon$, therefore
$\displaystyle 
x_{mk+j} \in 
%\left( f^j(\hat{x})-\frac{\varepsilon}{2}-\frac{\varepsilon}{2}, f^j(\hat{x})+\frac{\varepsilon}{2}+\frac{\varepsilon}{2} \right) = 
(f^j(\hat{x})-\varepsilon, f^j(\hat{x})+\varepsilon)$, $j=1, \dots, k$,
a.s.
\qed
%, which concludes the proof.
%\end{proof}
%%%%%%%%%%%%%%%%%%%%%%%%%

\section{Additive perturbations}
\label{sec:add}

In this section we investigate similar problems for stochastic equation with additive perturbations \eqref{4}, 
where $f$ satisfies Assumption \ref{as:slope}.
Our purpose remains the same: to achieve pseudo-stabilization of a blurred cycle $\{ f^j(x^*) \}$, $j=1, \dots, k$. 
Here $x^*$ is an arbitrary point $x^*\in \bigl(0, f(b)\bigr)$. 

Denoting again $g(x) := f^k(\nu x)$,
%\begin{equation}
%\label{h_def}
%g(x) := f^k(\nu x),
%\end{equation}
we can connect \eqref{4} to the equation with $x_{mk}=z_m$, $x_0=z_0>0$,
\begin{equation}
\label{4modified_1}
z_{m+1}= \max\left\{ g(z_m) + \ell\chi_{m+1}, 0\right\},\,   m\in \mathbb N.
\end{equation}
Let $x^* \in (0,f(b))$, $\nu= \Psi_k^{-1} (x^*)$, $\hat{x}= \nu x^* =f^{-k} (x^*) \in (0,b_k)$.
Note that $b_k/\nu>x^*$ and, for a fixed $\nu$, by Lemma~\ref{iterate}, 
$g(x)=f^k(\nu x)$ satisfies Assumption~\ref{as:slope}  
for $\nu x \in (0,b_k]$, so
$g(b_k/\nu)/(b_k/\nu) < g(x^*)/x^*$. Here the equality $g(x^*)=x^*$ is due to Lemma~\ref{lem:fixg}, Part 1. 
Thus $\displaystyle \frac{b_k}{\nu}-g\left( \frac{b_k}{\nu} \right) >0$. 
In addition, $g(x)>x$ for $x \in (0,x^*)$ and $g(x)>x$, $x \in (x^*, b_k/\nu)$.
For $\ell =0$, from monotonicity of $g$
on $(0,b_k/\nu)$, $x^*$ is an attractor of $g$ on $(0,b_k/\nu)$. Moreover, $g(x)<x$ for any $x>x^*$ implies $x^*$ is an 
attractor
for any $z_0 >0$. 
Our purpose is to choose $\ell>0$ small enough, to have $z_{m+1} \in (0,b_k/\nu)$, once $z_m$ is in this interval.

However, attractivity of a positive equilibrium in a deterministic case, in the presence of the zero equilibrium, does not imply that zero is a repeller in the stochastic case, see \cite{GDEA_2019,Medv} and references therein. 
Generally, with a positive probability, a solution can still stay in the right neighbourhood of zero.
Assumption~\ref{as:chibound} and its generalized version in Remark~\ref{rem:chibound} allow to make a conclusion on attractivity of $x^*$, a.s.

We choose $\delta_0>0$ satisfying
\begin{equation}
\label{add1}
\delta_0 < \min\left\{ \frac{b_k}{\nu}-g\left( \frac{b_k}{\nu} \right), \max_{x\in [0, x^{\ast}]} \left[ g(x)-x\right] \right\}.
\end{equation}

Define the numbers $y_1$, $y_2$, $\hat{x}_1$, $\hat{x}_2$ as
\begin{equation}
\label{add3}
\begin{split}
y_1 :=& \sup \left\{ \left. x \in [0, x^{\ast}] \right| g(x)-x \geq \delta_0 \right\}, %\in (0,  x^{\ast}), 
\\ \hat{x}_1:= & \nu y_1 \in (0,b_k),~y_1 \in (0,  x^{\ast}),\\
y_2 :=& \inf \left\{ \left. x \in [ x^{\ast},b_k/\nu] \right| g(x)-x \leq - \delta_0 \right\}, 
%\in ( x^{\ast},\frac{b_k}{\nu}),
\\ \hat{x}_2 := & \nu y_2 \in (0,b_k), ~y_2\in ( x^{\ast}, b_k/\nu). 
\end{split}
\end{equation}
According to the choice of $\delta_0$, the sets in \eqref{add3} are non-empty, so $y_1$,$y_2$, $\hat{x}_1$ and $\hat{x}_2$ are well defined.
Denote
\begin{equation}
\label{add4}
y_3 = \inf \left\{ \left. x \in [ x^{\ast},\infty) \right| g(x) - \delta_0 \leq y_1 \right\},
\end{equation}
where $y_3$ is assumed to be infinite if the set in the right-hand side of \eqref{add4} is empty.
As stated in \cite[Lemma 4.1]{DCDSB2017}, the numbers $y_1$, $y_2$ and $y_3$ defined by \eqref{add3} and 
\eqref{add4}, respectively, exist.

\begin{lemma} \cite[Theorem 4.5]{DCDSB2017}
\label{lemma:additive}
Let Assumptions~\ref{as:slope} and \ref{as:chibound}  hold,
$x^*\in (0, f(b))$  be an arbitrary point, $\nu=\nu(x^*)$ be chosen as in \eqref{def:nux*},
$g(x) = f^k(\nu x)$ and $\delta_0$ satisfy  \eqref{add1}.
Suppose that $y_1$, $y_2$, $y_3$ are denoted in  \eqref{add3} and \eqref{add4}, respectively,
and $z_m$ is a solution to equation \eqref{4modified_1} with an arbitrary $z_0>0$ and $\ell>0$ satisfying
%\begin{equation}
%\label{eq:l_min} 
$\ell \leq \delta_0$.
%\end{equation}
Then \\
\noindent(i) for each $\varepsilon_1>0$, there exists a random $\mathcal M(\omega)=\mathcal M(\omega, x_0, \ell, x^*, \varepsilon_1)$ such that
for $m \ge \mathcal M(\omega)$ we have, a.s. on $\Omega$,
\begin{equation}
\label{add50a0}
y_1\le z_m\le  y_2+\varepsilon_1;
\end{equation}
\noindent(ii) for each $\varepsilon_1>0$ and $\gamma \in (0, 1)$, there is a nonrandom 
number  $M=M(\gamma, x_0,\ell,x^*,\varepsilon_1)$
such that
\begin{equation}
\label{add50a}
\mathbb P \{y_1\le z_m\le  y_2+\varepsilon_1, \,\, \text{\rm{for}} \,\,  m\ge M\}>\gamma;
\end{equation}
(iii)
we have
%\begin{equation}
%\label{add5a}
$\displaystyle \liminf_{m \to \infty} z_m \geq y_1$, $\displaystyle \limsup_{m \to \infty} z_m \leq y_2$,  
\text{a.s}.
%\end{equation}
%\end{enumerate}
\end{lemma}

Another result that will be used in future is also stated below.
It illustrates that a solution will eventually be in any arbitrarily small neighborhood 
of %the point 
$x^{\ast}$ with an arbitrarily close to one probability and 
will further be used in the proof of Theorem~\ref{theorem:l0add}.

\begin{lemma}\cite[Theorem 4.6]{DCDSB2017}
\label{prop:l0add}
Let Assumptions~\ref{as:slope} and \ref{as:chibound}  hold,  $z_0>0$ be an arbitrary initial value, 
$x^*\in (0, f(b))$  be an arbitrary point, $\nu=\nu(x^*)$ be chosen as in \eqref{def:nux*}.
Then, for each $\varepsilon>0$ and $\gamma\in (0, 1)$, we can find $\delta_0$ such that for the solution $z_m$ to 
\eqref{4modified_1} with  $\ell\le \delta_0$,   and for  some
nonrandom  $M=M(\gamma, x_0, \ell, x^*, \varepsilon) \in \mathbb N$, we have
\\
$\mathbb P\{z_m\in (x^*-\varepsilon, x^*+\varepsilon) \,\,  \forall m \ge M \}\ge \gamma$.
\end{lemma}

%!!!!!!!!!!!!!!!!!!!!!!!!!!!!!!!!!!!!!!!!

This leads to two main results %on stable blurred $k$-cycles of 
for \eqref{4},
Lemma~\ref{lemma:additive} implying Theorem~\ref{th:additive} on a.s. convergence to a blurred cycle, and Lemma~\ref{prop:l0add} yielding 
Theorem~\ref{theorem:l0add} on the convergence with a prescribed close to one probability.

%!!!!!!!!!!!!!!!!!!!!!!!!!!!!!!!!!!!!!!!!!!!!!!!!

\begin{theorem}
\label{th:additive}
Let Assumptions~\ref{as:slope} and \ref{as:chibound}  hold, 
$\hat{x} \in  (0, b_k)$  be an arbitrary point, $x^*=f^k(\hat{x})$, $\nu=\nu(x^*)$ be chosen as in \eqref{def:nux*},
%$h$ be defined as in \eqref{h_def}, 
and $\delta_0$ satisfy  \eqref{add1}.
Suppose that $\hat{x}_1$ and  $\hat{x}_2$ are defined as in  \eqref{add3}, 
and $x_n$ is a solution to %equation 
\eqref{4} with an arbitrary $x_0>0$ and $\ell>0$ satisfying 
$\ell \leq \delta_0$.
Then 
%\begin{enumerate}
%\item [(i)] 
\\
(i) For any $\varepsilon>0$, there exists a random ${\mathcal M}(\omega)={\mathcal M}(\omega, x_0, \ell, \hat{x}, \varepsilon)$ such that
for $m \ge \mathcal M(\omega)$ we have, a.s. on $\Omega$, $\displaystyle f^j(\hat{x}_1) \leq x_{km+j} \le  f^j 
(\hat{x}_2)+\varepsilon$, $j=1, \dots k$.
%\begin{equation}
%\label{add_for_4}
%f^j(\hat{x}_1) \le x_{km+j} \le  f^j (\hat{x}_2)+\varepsilon, \quad j=1, \dots k;
%\end{equation}
%\noindent(ii) 
%\item [(ii)] 
\\
(ii)
For each $\varepsilon>0$ and $\gamma \in (0, 1)$, there is a nonrandom 
number  $M=M(\gamma, x_0,\ell,\hat{x},\varepsilon)$
such that, for $j=1, \dots k$,
\begin{equation}
\label{add50afor_x} 
\mathbb P \{f^j(\hat{x}_1) \le x_{km+j} \le  f^j (\hat{x}_2)+\varepsilon, \, m\ge M\}>\gamma, 
\end{equation}
\\
(iii)
%\item [(iii)] 
For a solution $x_n$ of \eqref{4} we have, a.s., for $j=1, \dots k$,
\begin{equation}
\label{add5afor_x}
\vspace{-2mm}
%\begin{array}{l}
\displaystyle
\liminf_{n \to \infty} x_{km+j} \geq f^j (\hat{x}_1),  
\limsup_{n \to \infty} x_{km+j} \leq f^j(\hat{x}_2).
%\quad j=1, \dots k, \quad 
%\text{~a.s}, j=1, \dots k.
%\end{array}
\end{equation}
%$j=1, \dots k$.
%\end{enumerate}
\end{theorem}
{\bf Proof.}
Recall from \eqref{add3} that $\nu y_1= \hat{x}_1$, $\nu y_2= \hat{x}_2$. 
From continuity and monotonicity of $f$, for any $\varepsilon>0$, 
there is a $\varepsilon_1>0$ such that  \eqref{add50a0} implies 
%for a random $\mathcal M(\omega)=\mathcal M(\omega, x_0, \ell, \hat x, \varepsilon_1)$  
%for $m \ge \mathcal M(\omega)$,
\begin{equation}
\label{add_1star}
f^j(\hat{x}_1) \leq f^j(\nu z_m) \leq f^j(\hat{x}_2)+ \varepsilon, ~~j=1, \dots, k.
\end{equation}
We have
\begin{equation}
\label{add_2star}
x_{mk}=z_m, \quad x_{mk+j}=f^j(\nu z_m), ~~ j=1, \dots, k.
\end{equation}
(i) Choosing this $\varepsilon_1$ as in (i) of Lemma~\ref{lemma:additive}, 
we find $\mathcal M(\omega)=\mathcal M(\omega, x_0, \ell, x^*, \varepsilon_1)$ such that \eqref{add50a0}, and  thus \eqref{add50afor_x} are 
satisfied. 
\\
(ii) Further, (ii) in Lemma~\ref{lemma:additive} implies
for $M=M(\gamma, x_0,\ell,x^*,\varepsilon_1)$ inequality
\eqref{add50a}. Thus by \eqref{add_1star} and \eqref{add_2star} we have
$
\mathbb P \{f^j(\hat{x}_1) \le x_{km+j} \le  f^j(\hat{x}_2)+\varepsilon, \,\, \text{\rm{for}} \,\,  m\ge M\}
$ \\ $
\geq {\mathbb P} \{\hat{x}_1\le z_m\le  y_2+\varepsilon, \,\, \text{\rm{for}} \,\,  m\ge M\}>\gamma.
$
\\
(iii) As $x_{mk+j}$ and $z_m$ are connected with \eqref{add_2star}, 
application of Part (iii)  in Lemma~\ref{lemma:additive} immediately implies \eqref{add5afor_x}. 
\qed
% and concludes the proof.
%\end{proof}
%\medskip

\begin{theorem}
\label{theorem:l0add}
Let Assumptions~\ref{as:slope} and \ref{as:chibound}  hold,  $x_0>0$, $\hat{x}\in (0, b_k)$  be an arbitrary point, $x^*=f^k(\hat{x})$, 
$\nu=\nu(x^*)$ be chosen as in \eqref{def:nux*}.
Then, for each $\varepsilon>0$ and $\gamma\in (0, 1)$, we can find $\delta_0$ such that for the solution $x_n$ to 
\eqref{4} with  $\ell\le \delta_0$,   and for  some
nonrandom  $M=M(\gamma, x_0, \ell, \hat x, \varepsilon) \in \mathbb N$, $j=1, \dots, k$, we have
\\
%\begin{equation}
%\label{add_3star}
%\vspace{-0.3cm}
${\mathbb P} \left\{ x_{km+j} \in \left( f^j(\hat{x})-\varepsilon, f^j(\hat{x})+\varepsilon \right) \,  
%\text{\rm for all } 
\forall
m \ge M \right\}\ge \gamma$.
%\end{equation}
\end{theorem}
{\bf Proof.} Let us choose $\varepsilon_1$ such that \eqref{add_1star} is satisfied, fix $\gamma \in (0,1)$ and find
$M=M(\gamma, x_0, \ell, \hat x, \varepsilon_1) \in \mathbb N$ as in Lemma~\ref{prop:l0add}. Then
$\displaystyle
\mathbb P\{\nu z_m\in (\hat{x}-\varepsilon_1, \hat{x}+\varepsilon_1) \,\,  \text{\rm for all $m \ge M$}\}\ge \gamma$,
which by \eqref{add_2star} implies the statement of the theorem. %\eqref{add_3star}. 
\qed

%This concludes the proof.
%\medskip

\section{Examples}
\label{sec:ex}

%Before proceeding  to examples consider equation where we  combine multiplicative and additive noise 
%\begin{equation}
%\label{5_new}
%x_{n+1}= \left\{ \begin{array}{ll} \displaystyle f\left( (\nu + \ell_1\chi_{m+1})x_n \right) + \ell_2\chi_{m+1},
%& n=mk, ~m\in {\mathbb N},\\
%f(x_n), & n \neq mk,  \end{array} \right.
%\end{equation}
%$$x_0>0, \quad n\in \mathbb N, \quad \nu \in (0,1], \quad k\in {\mathbb N}.$$

We consider \eqref{5_new} combining multiplicative and additive noise.
Similarly to the previous theorems, the following more general result can be obtained. However, the  proof is 
long and technical and does not include any new ideas. Therefore we do not present it, but only 
illustrate stated below Proposition~\ref{prop:combined} with computer simulations.

\begin{proposition}
\label{prop:combined}
Let Assumptions~\ref{as:slope} and \ref{as:chibound}  hold,  $x_0>0$, $\hat{x} \in (0, b_k)$  be an arbitrary point,
$x^*=f^k(\hat{x})$,
$\nu=\nu(x^*)$ be chosen as in \eqref{def:nux*}.
Then, for each $\varepsilon>0$ and $\gamma\in (0, 1)$, we can find $\delta_1$ and $\delta_2$ such that for the solution $x_n$ to
\eqref{5_new} with  $\ell_1 \le \delta_1$, $\ell_2 \le \delta_2$  and for  some   
nonrandom  $M=M(\gamma, x_0, \ell_1, \ell_2, \hat x, \varepsilon) \in \mathbb N$, we have
%\begin{equation*}
%\label{add_4star}
\\
$\displaystyle {\mathbb P} \left\{ x_{km+j} \in \left( f^j(\hat{x})-\varepsilon, f^j (\hat{x})+\varepsilon \right)
\forall  m \ge M \right\}\ge \gamma$,
$j=1, \dots, k$.
%\end{equation*}
\end{proposition}

Now we present examples of application of noisy  PF control method  to create a stable equilibrium or stable 
$k$-cycle in the neighborhood of nonzero point $K$. In all case  noises $\chi $ are 
continuous uniformly distributed on $[-1,1]$. In all the simulations five runs with the same initial value are illustrated,
with $n$ on the $x$-axis and $x_n$ (for all the five runs) on $y$-axis.

\begin{example}
\label{ex:1}
Let us apply PF control to  the Ricker model \eqref{eq:ricker}. For $r=2.8$, the non-controlled map is chaotic. We consider $\nu=0.002$,  noise applied every third step.
%First, we consider \eqref{3} with 
%$\nu=0.03$, $\ell=0.001$ and noise applied every second step, 
%and then 
%$\nu=0.002$, $\ell_1=0.0001$, noise applied every third step.
For \eqref{3} with  $\ell_1=0.0001$ we observe a blurred stable 
3-cycle, % in latter case, see Figs.~\ref{figure1a}
see Fig.~\ref{figure2a}, left.
Next, we simulate additive noise as in \eqref{4}.
We observe a blurred stable %2-cycle  and a blurred 
3-cycle %, respectively, 
with similar amplitudes for larger $\ell_2$, see 
Fig.~%~\ref{figure1a} and 
\ref{figure2a}, right. For the combined noise as in \eqref{5_new}, the results of the runs are similar to Fig.~\ref{figure2a}, left.

\begin{figure}[ht]
\centering
%!!!!!!!!!!!!!!!!!!!!!!!!!!!!!!!!!!!!!!!!!!!!
\includegraphics[height=.125\textheight]{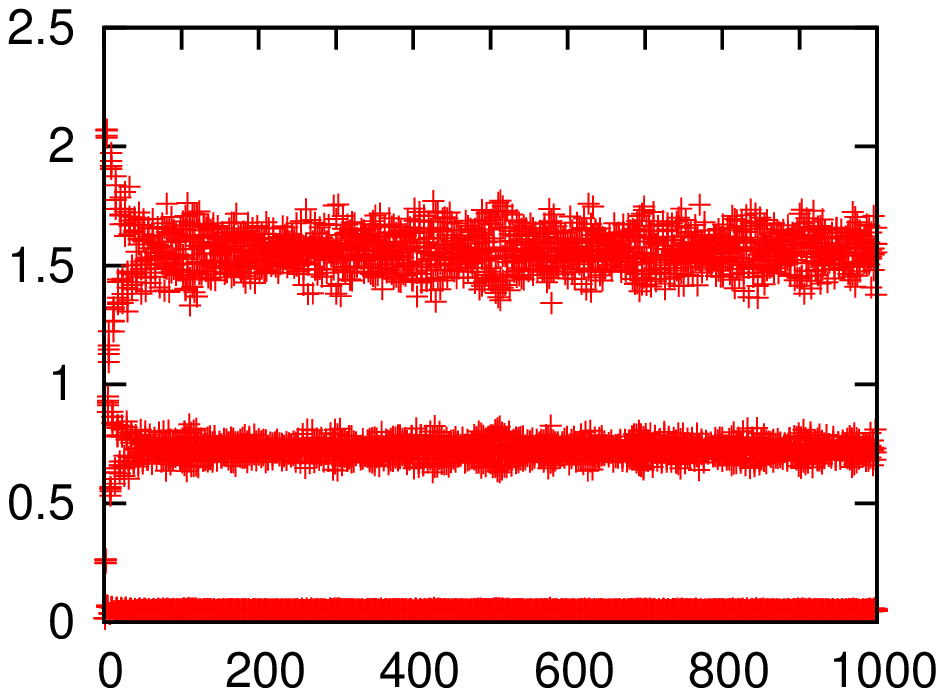}
%???????????????????????????????????
%\hspace{6mm}
%!!!!!!!!!!!!!!!!!!!!!!!!!!!!!!!!!!!!!!!!!!!!
\includegraphics[height=.125\textheight]{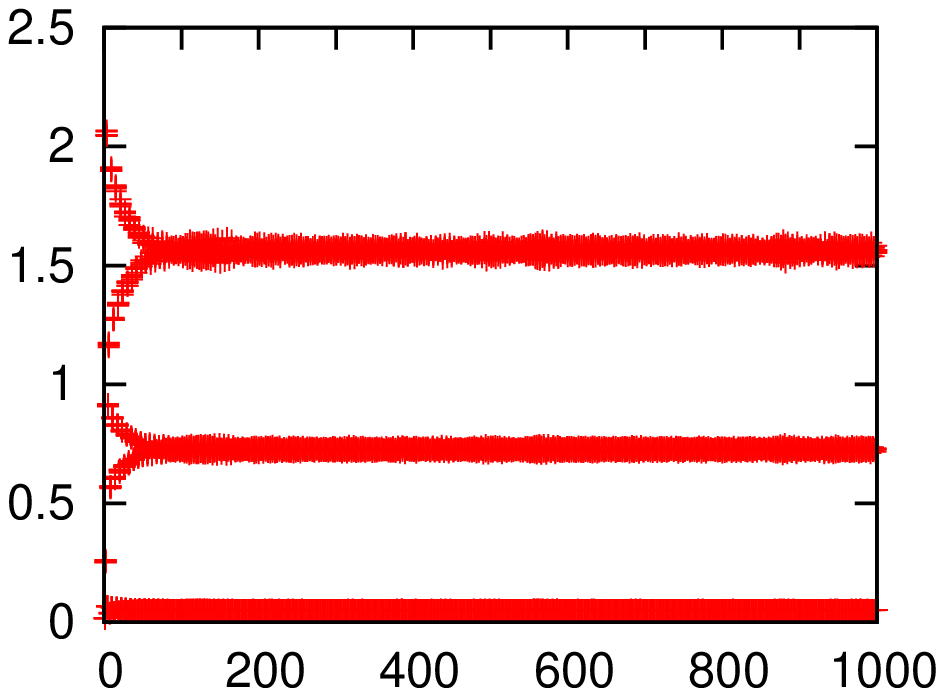}
%???????????????????????????????????
%\hspace{6mm} 
%\includegraphics[height=.08\textheight]{Prog/Impulse/Rick_both_n_3_l_1_0_0001_l_2_0_0005_r_2_8_nu_0_002.ps}
\caption{Solutions of the Ricker difference equation with $f=f_1$ as in
(\protect{\ref{eq:ricker}}), $r=2.8$, $x_0=0.5$,
$k=3$, $\nu=0.002$, $n=0, \dots, 1000$
and (left)  (\protect{\ref{3}}) with $\ell_1=0.0001$, 
%(middle) equation (\protect{\ref{4}})
(right)  (\protect{\ref{4}})
with  $\ell_2=0.0005$.
%(right) equation (\protect{\ref{5_new}}) with $\ell_1=0.0001$, $\ell_2=0.0005$.
%Everywhere $x_0=0.5$.
}
\label{figure2a} 
\end{figure}

\end{example}

\begin{example}
\label{ex:2}
Consider a particular case of
\eqref{eq:milton}, see \cite{Liz_CAMWA,DCDSB2017},
\begin{equation}
\label{eq:CAMWA}
f(x) = x \left( 0.55 + \frac{3.45}{1+x^9} \right), \quad x\ge0.
\end{equation}
We apply PF with %$k=2$ and 
$k=3$ to the three cases: the multiplicative noise, as in \eqref{3},
the additive noise, as in \eqref{4}, and the combined noise as in \eqref{5_new},
see Fig.~\ref{figure4a}.

\begin{figure}[ht]
\centering
%!!!!!!!!!!!!!!!!!!!!!!!!!!!!!!!!!!!!!!!!!!!!
\includegraphics[height=.12\textheight]{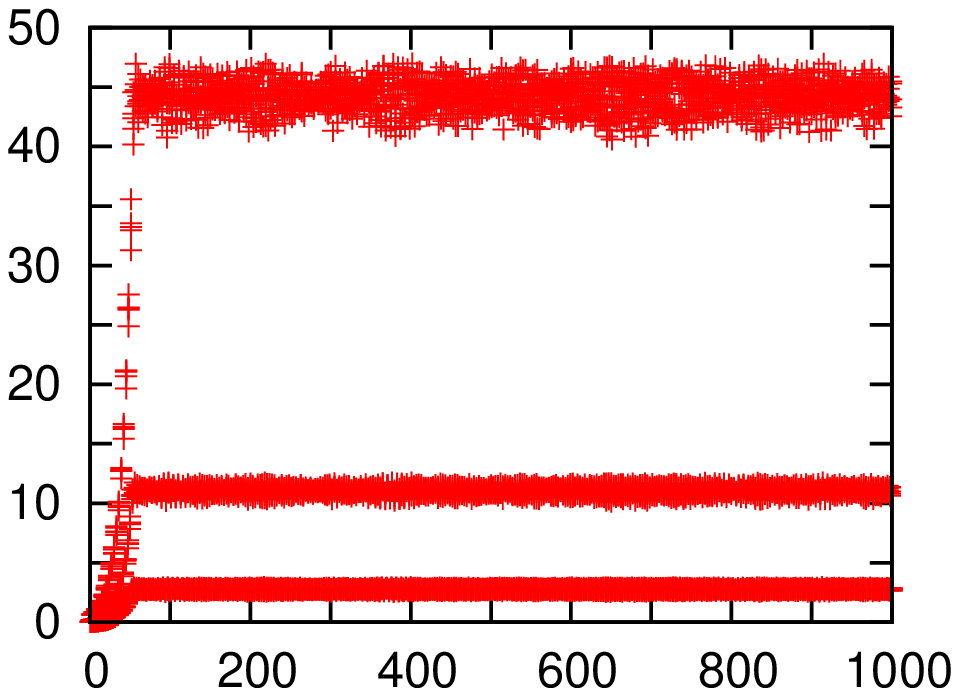}
%???????????????????????????????????
%\hspace{6mm}
%!!!!!!!!!!!!!!!!!!!!!!!!!!!!!!!!!!!!!!!!!!!!
\includegraphics[height=.12\textheight]{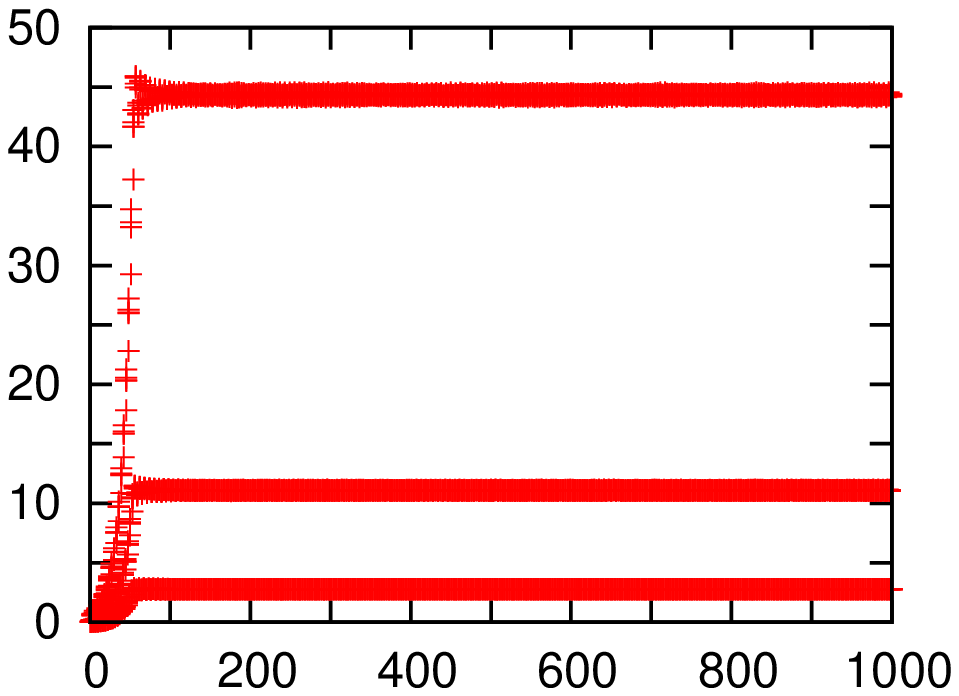}
%???????????????????????????????????
%\hspace{6mm}   
\includegraphics[height=.12\textheight]{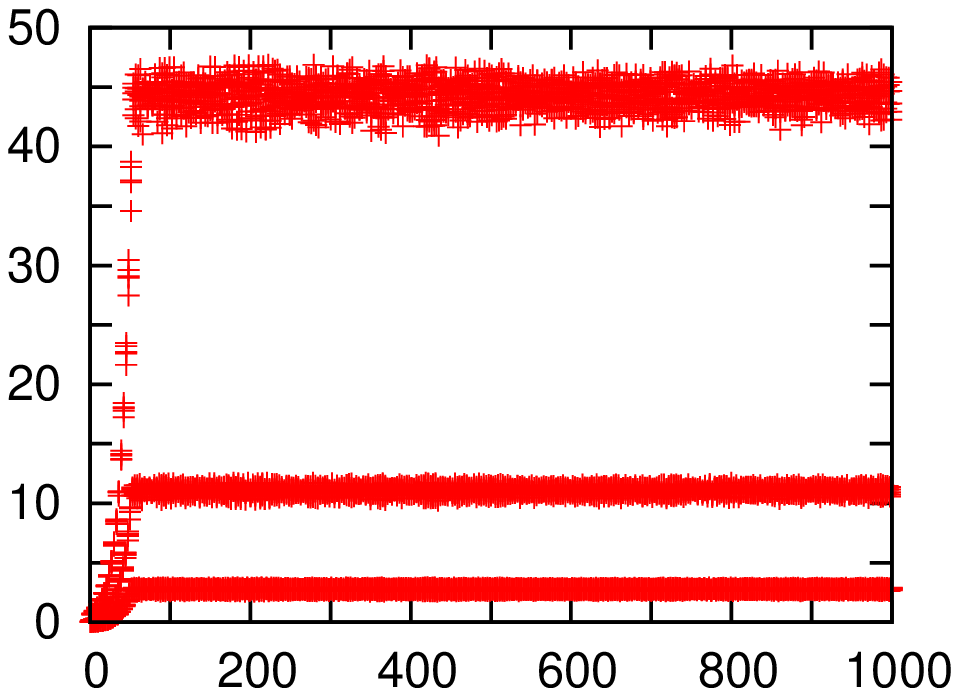}
\caption{Solutions $x_n$ vs. $n$  with $f$ as in
(\protect{\ref{eq:CAMWA}}),
$k=3$, $\nu=0.02$, $x_0=0.5$,
and (top left)  (\protect{\ref{3}}) with $\ell=0.0005$, (top right)  (\protect{\ref{4}})
with  $\ell=0.005$,
(bottom)  (\protect{\ref{5_new}}) with $\ell_1=0.0005$ and $\ell_2=0.005$.
%Everywhere $x_0=0.5$.
}
\label{figure4a}
\end{figure}
\end{example}

The standard PF control moves a positive equilibrium towards zero; applied at every $k$th step, it leads to a stable 
cycle
in a right neighbourhood of zero. Now we modify this method choosing a positive equilibrium $K_1$ instead of zero.  
We apply PF control method  to create a stable equilibrium or $k$-cycle in the nighbourhood of nonzero point $K_1$. 
%Different transformations are used to the left and to the right sides of  $K_1$, 
%each $1/\nu$- constructs  the solution to $K_1$.
The non-shifted PF control brings the state variable $1/\nu$ times closer to zero. We mimic this idea for a 
shifted version assuming that the state variable is proportionally moved to the fixed $K_1$.
The controlled equation has the form %equilibrium $K_1$:
$x_{n+1} =  f(K_1+\nu (x_n-K_1))-K_1+K_1 = f( \nu x_n +(1-\nu)K_1)$, $x_n\geq K_1$, 
$x_{n+1} =  K_1-\left[K_1-f \left( K_1-\nu(K_1-x_n) \right) \right] = f( \nu x_n +(1-\nu)K_1)$,
$x_n\in (0,K_1)$. Thus
\begin{equation}
 \label{eq:TOC}
%\begin{split}
x_{n+1} =  %f(K_1+\nu (x_n-K_1))-K_1+K_1 = 
f( \nu x_n +(1-\nu)K_1).%, ~~x_n\geq K_1, \\
%x_{n+1} =  K_1-\left[K_1-f \left( K_1-\nu(K_1-x_n) \right) \right] = f( \nu x_n +(1-\nu)K_1), ~~x_n\in (0,K_1).
%\end{split}
\end{equation}

\begin{example}
\label{ex:3}
Define
 \begin{equation}
 \label{def:efx}
f(x):=\frac 92 x^2(1-x), \quad x\in [0, 1].
\end{equation}
The maximum value of  $f_{\max}$ is achieved at  $x_{\max}=\frac{2}{3}$, $f(x_{\max})=\frac{2}{3}$, the inflection point is $x^{\ast}=\frac 13$, 
$f''(x)>0$ for $x\in (0, \frac 13)$ and $f''(x)<0$ for $x\in (\frac 13, 1)$,   $f$ has two positive equilibrium points 
$K_1=\frac 13$,  $K_2=\frac 23$ and  $f'\left(\frac 13\right)=\frac 32>1$.

Consider  a modification of PF method ``centered'' at  $K_1=1/3$, see \eqref{eq:TOC}.
%Applying the results of Sections 3 and 4 
%after the  transformations $x=z+K_1$  and $x=K_1-y$,  we arrive at the equation
%\begin{equation}
%\label{other_centre}
%x_{n+1}=f(\nu x_n+(1-\nu) K_1), \quad x_0\in \left(0,   2/3\right).
%\end{equation} 
It can be shown that, for $\nu\in (2/3, 1)$, equation  \eqref{eq:TOC} has  two positive locally  stable 
equilibrium points on both sides of  $K_1$,  each attracts a solution  $x_n$ with corresponding position of 
$x_0$ around $K_1$, see bifurcation diagram on Fig~\ref {figure5a}.

%set xlabel "c"
%set output "bifurc_general.ps"
%plot "bifurc_general", 1/3.0
%set output "bifurc_low.ps"
%plot "bifurc_low", 1/3.0
%set output "bifurc_high.ps"
%plot "bifurc_high", 1/3.0
%set output "bifurc_general_coef_6.ps"
%plot "bifurc_general_coef_6", 0.5*(1-sqrt(1.0/3.0))

%set xlabel "c"
%set output "bif_gen.ps"
%plot [0:0.9] [0:1] "bif_gen", 1/3.0
%set output "bif_low.ps"
%plot [0:0.9] [0:1] "bif_low", 1/3.0
%set output "bif_high.ps"
%plot [0:0.9] [0:1] "bif_high", 1/3.0
%set output "bif_gen_coef_6.ps"
%plot [0:0.9] [0:1] "bif_gen_coef_6", 0.5*(1-sqrt(1.0/3.0))

\begin{figure}[ht]
\centering
%!!!!!!!!!!!!!!!!!!!!!!!!!!!!!!!!!!!!!
\includegraphics[height=.16\textheight]{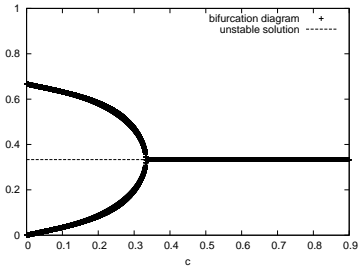}
%???????????????????????????????????????
%\includegraphics[height=.16\textheight]{Prog/Impulse/bif_low.ps}
%\includegraphics[height=.16\textheight]{Prog/Impulse/bif_high.ps}
%\includegraphics[height=.14\textheight]{Prog/Impulse/bifurc_general.ps}
%\includegraphics[height=.14\textheight]{Prog/Impulse/bifurc_low.ps}
%\includegraphics[height=.14\textheight]{Prog/Impulse/bifurc_high.ps}
\caption{Bifurcation diagram for (\protect{\ref{eq:TOC}}) with $f$ as in (\protect{\ref{def:efx}}), $c=1-\nu$ changing from zero to 0.9  and $x_0$ changing from 0 to 1.
We get an upper branch if  $x_0$ changes from 1/3 to 1 and the lower branch if it changes from zero to 1/3.
% (left), 
%from zero to 1/3 (middle) and from 1/3 to 1 (right).}
}
\label{figure5a}
\end{figure}

Note that \eqref{eq:TOC} is a particular case of Target Oriented Control
\cite{Dattani}, sufficient conditions for stabilization of $K_1$ in \eqref{eq:TOC} were obtained in
\cite{TPC}. A modification of PF method is responsible for the left part of the diagram (bistability) while 
\cite{TPC} gives an exact bound $c^*$ such that for $c\in (c^*,1)$, all solutions of \eqref{eq:TOC} with $\nu:=1-c$ and  $x_0 \in (0,1)$ converge to $K_1=1/3$.

We introduce multiplicative noise in \eqref{eq:TOC} to get for any $k\in {\mathbb N}$, $\nu \in (0,1]$,
\begin{equation}
\label{eq5_7}
x_{n+1}= \left\{ \begin{array}{l} \displaystyle %f\left( (\nu + \ell_1\chi_{m+1})x_n \right),
f((\nu+\ell_1\chi_{m+1}) x_n \\ +(1-\nu-\ell_1\chi_{m+1}) K_1),
~ n \mid k, \\ %~m,n \in {\mathbb N}_0,~ x_0>0,\\
f(x_n),  n \not\,\mid k, \\ n \in {\mathbb N}_0,~ x_0>0.  \end{array} \right.
%, ~ \nu \in (0,1],~ k\in {\mathbb N}.
\end{equation}
A multiplicative noise with small $\ell_1$ does not change this type of behavior, as illustrated in 
Fig.~\ref{figure6a}. 
This also holds when  coefficient $\ell_2$ of the additive noise  is relatively small and $x_0$ is relatively far from 
$K_1$,  
see Fig~\ref {figure7}, left and middle.  However, when  $\ell_2$  increases (in some limits), the solution started on the left of 
$K_1$ and close enough to $K_1$, is attracted to both equilibrium solutions, on the left and on the right of $K_1$, see Fig.~\ref{figure7}, right. 
The same holds when  $x_0>K_1$.
Fig.~ %\ref{figure8} and  
\ref{figure9} 
illustrates  construction of stable %two and 
three-cycles 
when the initial value is taken on both sides of $K_1$.

%set term postscript portrait
%set size 1.0,0.5
%set output "Allee_mult_nu_0_7_l_0_0005_x_0_35.ps"
%plot "Allee_mult_nu_0_7_l_0_0005_x_0_35"         
%set output "Allee_mult_nu_0_7_l_0_0005_x_0_6.ps"
%plot "Allee_mult_nu_0_7_l_0_0005_x_0_6"
%set output "Allee_mult_nu_0_7_l_0_0005_x_0_2.ps"
%plot "Allee_mult_nu_0_7_l_0_0005_x_0_2"

%set output "Allee_add_nu_0_7_l_0_001_x_0_4.ps"
%plot "Allee_add_nu_0_7_l_0_001_x_0_4"
%set output "Allee_add_nu_0_7_l_0_001_x_0_3.ps"
%plot "Allee_add_nu_0_7_l_0_001_x_0_3"
%set output "Allee_add_nu_0_8_l_0_01_x_0_33.ps"
%plot "Allee_add_nu_0_8_l_0_01_x_0_33"

%set term postscript portrait
%set size 0.6,0.3
%set output "Allee_mult_nu_0_7_l_0_0005_x_0_35_new.ps"
%plot "Allee_mult_nu_0_7_l_0_0005_x_0_35" notitle
%set output "Allee_mult_nu_0_7_l_0_0005_x_0_2_new.ps"
%plot "Allee_mult_nu_0_7_l_0_0005_x_0_2" notitle
%set output "Allee_add_nu_0_7_l_0_001_x_0_4_new.ps"
%plot "Allee_add_nu_0_7_l_0_001_x_0_4" notitle
%set output "Allee_add_nu_0_7_l_0_001_x_0_3_new.ps"
%plot "Allee_add_nu_0_7_l_0_001_x_0_3" notitle
%set output "Allee_add_nu_0_8_l_0_01_x_0_33_new.ps"
%plot "Allee_add_nu_0_8_l_0_01_x_0_33" notitle

\begin{figure}[ht]
\centering
%!!!!!!!!!!!!!!!!!
\includegraphics[height=.12\textheight]{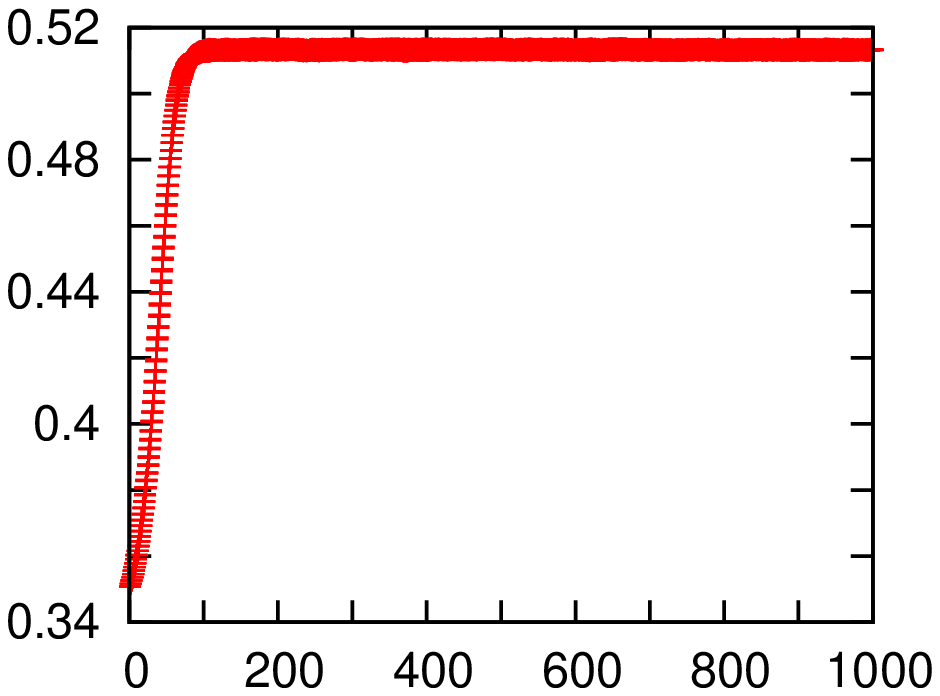}
%??????????????????????????
%\hspace{6mm}
%\includegraphics[height=.08\textheight]{Prog/Impulse/Allee_mult_nu_0_7_l_0_0005_x_0_6.ps}
%\hspace{6mm}
%!!!!!!!!!!!!!!!!!!!!!!!
\includegraphics[height=.12\textheight]{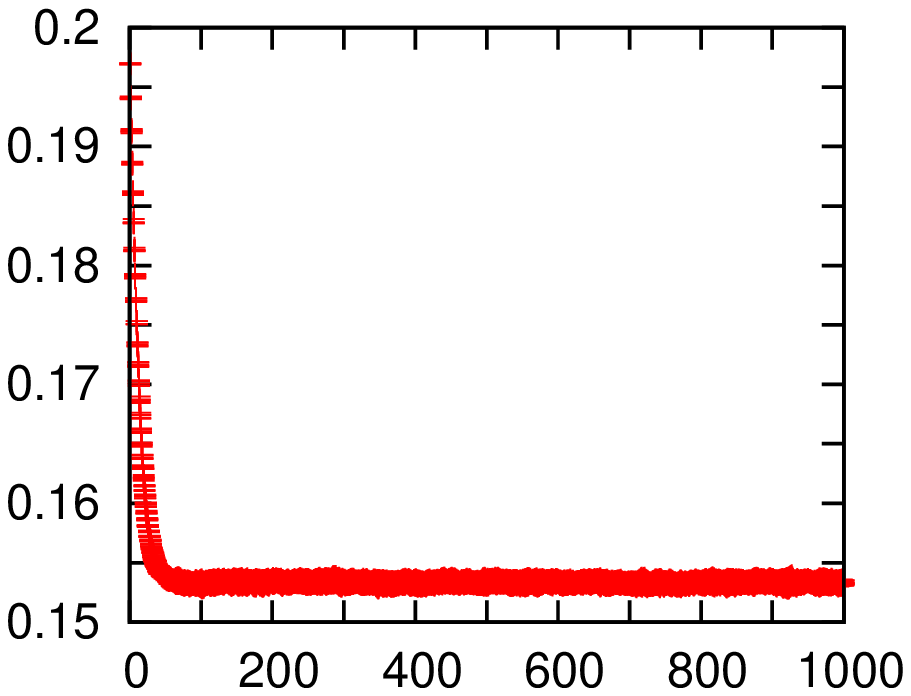}
%??????????????????????????
\caption{Five runs of the difference equation with $f$ as in (\protect{\ref{def:efx}}), 
%as in (\protect{\ref{other_centre}}), 
multiplicative noise with $\ell=0.0005$, $\nu=0.7$ and
(left) $x_0=0.35$, (middle) $x_0=0.6$, (right) $x_0=0.2$. 
}
\label{figure6a}
\end{figure}

\begin{figure}[ht]
\centering
%!!!!!!!!!!!!!!!!!!!!!!!!!!!!!!!!!!!!!!!!!!!!
\includegraphics[height=.12\textheight]{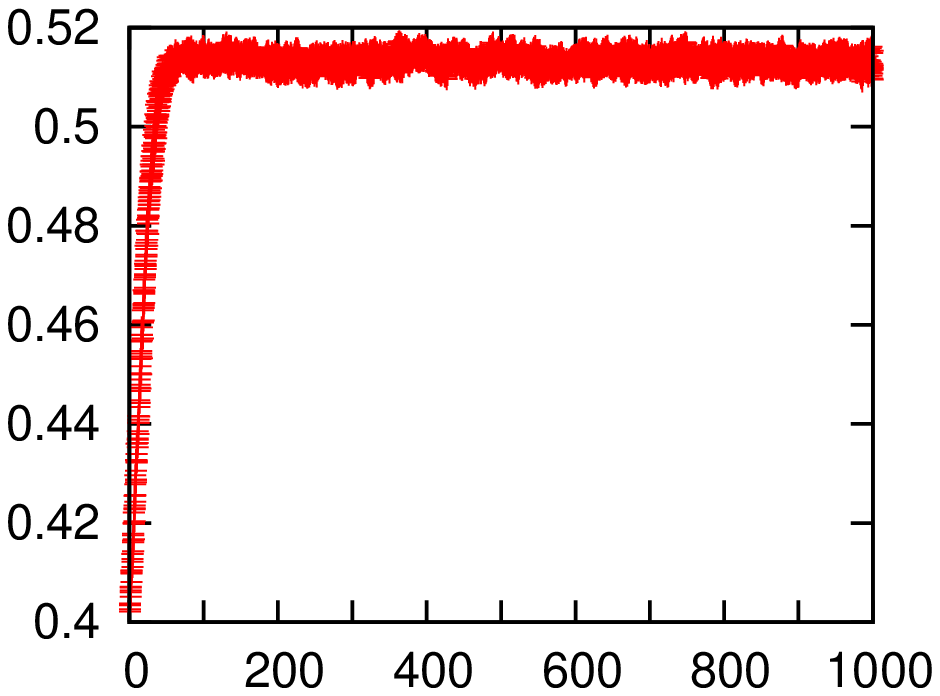}
%???????????????????????????????????
%\hspace{6mm}
%!!!!!!!!!!!!!!!!!!!!!!!!!!!!!!!!!!!!!!!!!!!!
\includegraphics[height=.12\textheight]{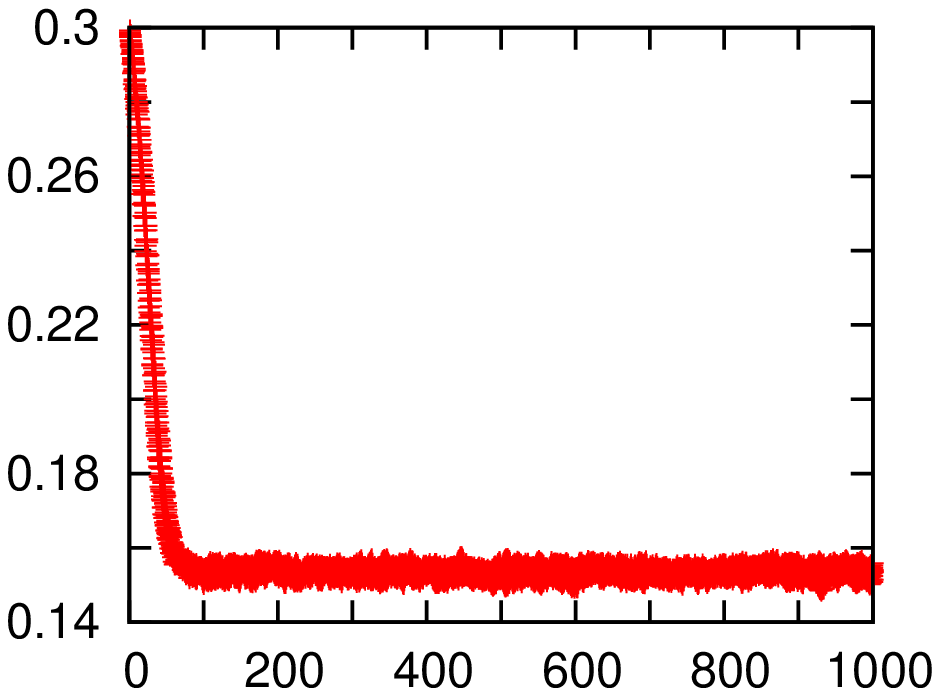}
%???????????????????????????????????
%\hspace{6mm}
%!!!!!!!!!!!!!!!!!!!!!!!!!!!!!!!!!!!!!!!!!!!!
\includegraphics[height=.12\textheight]{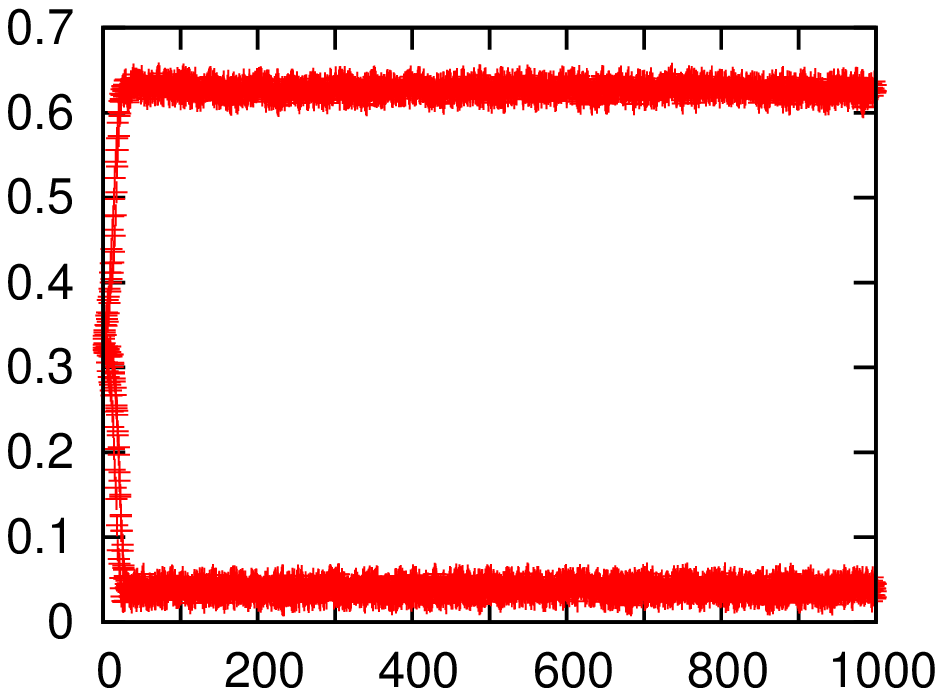}
%???????????????????????????????????
\caption{For difference equation (\protect{\ref{eq5_7}}) with $f$ as in (\protect{\ref{def:efx}}), with additive noise (left) $\nu=0.7$, $\ell=0.001$, $x_0=0.4$,
(middle) $\nu=0.7$, $\ell=0.001$, $x_0=0.3$, (right) $\nu=0.8$, $\ell=0.01$, $x_0=0.33$.
}
\label{figure7}
\end{figure}

\begin{figure}[ht]
\centering
%!!!!!!!!!!!!!!!!!!!!!!!!!!!!!!!!!!!!!!!!!!!!
\includegraphics[height=.12\textheight]{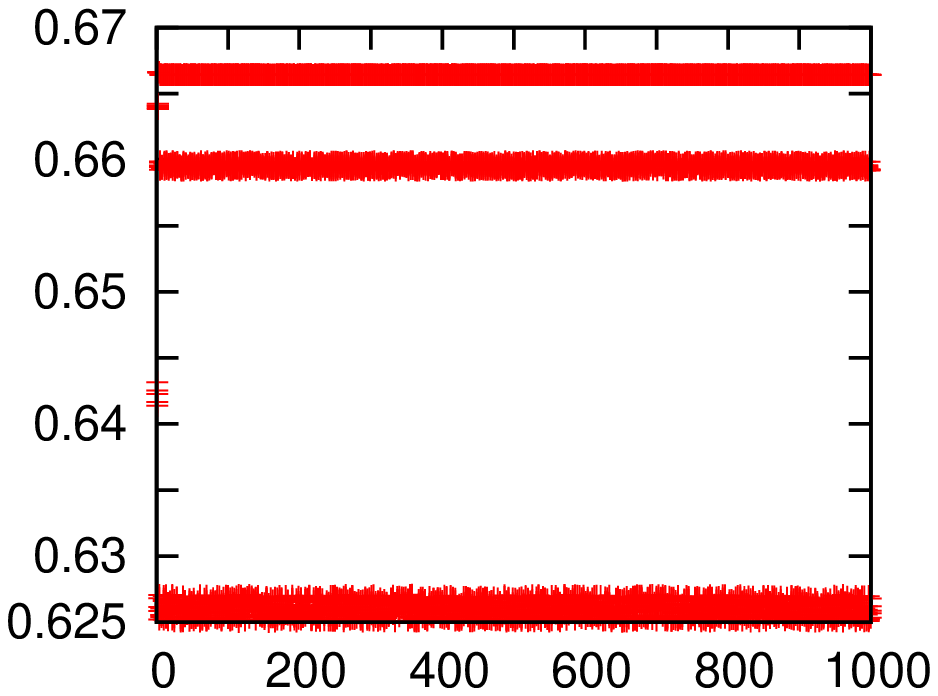}
%???????????????????????????????????
~~%\hspace{6mm}
%!!!!!!!!!!!!!!!!!!!!!!!!!!!!!!!!!!!!!!!!!!!!
\includegraphics[height=.12\textheight]{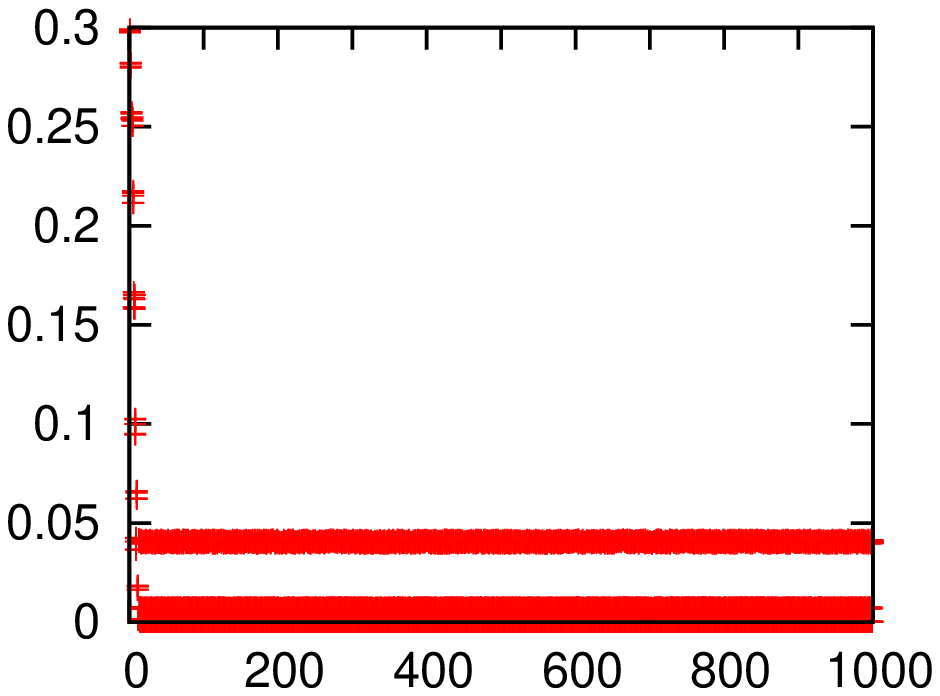}
%???????????????????????????????????
\caption{Five runs of difference equation (\protect{\ref{eq5_7}}) with $f$ as in (\protect{\ref{def:efx}}), PF control applied every third step,
%as in (\protect{\ref{other_centre}}),
multiplicative noise with $\ell_1=0.0001$, additive noise with $\ell_2=0.001$, $\nu=0.7$ and
(left) $x_0=0.3$, (right) $x_0=0.7$.
}
\label{figure9}
\end{figure}
\end{example}

%%%% For Automatica, omit ??? part of ??? the rest

\begin{example}
\label{ex:3a}
Define now  
 \begin{equation}
 \label{def:efx6}
f(x):=6 x^2(1-x), \quad x\in [0, 1],
\end{equation}
which has a  positive equilibrium $K_1\approx 0.211 <1/3$. Note that for $f$ as in \eqref{def:efx6},  
the results of Sections 3-4 can be applied for $x_n$ to the left of $K_1$, see the  bifurcation diagram in Figure \ref{figure9}.

\begin{figure}[ht]
\centering
%!!!!!!!!!!!!!!!!!!!!!!!!!!!!!!!!!!!!!!!!!!!!!!!!
\includegraphics[height=.16\textheight]{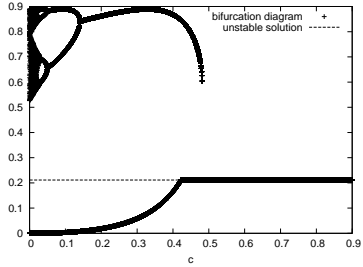}
%????????????????????????????????????????????????
\caption{Bifurcation diagram with $1-\nu$ changing from zero to 0.9 for the map $f$ as in (\protect{\ref{def:efx6}}).}
\label{figure9a}
\end{figure}

%Fig.~\ref{figure10}, low left,  demonstrates that a control with a  small multiplicative noise  keeps the solution on the 
%left of $K_1$, but when there is also an additive noise, a solution is attracted to both equlibriums,  
%Fig.~\ref{figure10}, upper left.  Even though $f$  does not satisfy our assumptions  on the interval right of $K_1$, 
%both pictures on  Fig \ref{figure10}, right,  demonstrates attraction to the equilibrium.

Fig.~\ref{figure11} illustrates a construction of a stable 2-cycle with multiplicative and additive noise. 
The left-side 
pictures, where the  initial value $x_0<K_1$, show  a 2-cycle, while the right-side pictures, where  $x_0>K_1$, produce 
a 3-cycle.

\begin{figure}[ht]
\centering
%!!!!!!!!!!!!!!!!!!!!!!!!!!!!!!!!!!!!!!!!!!!!
\includegraphics[height=.11\textheight]{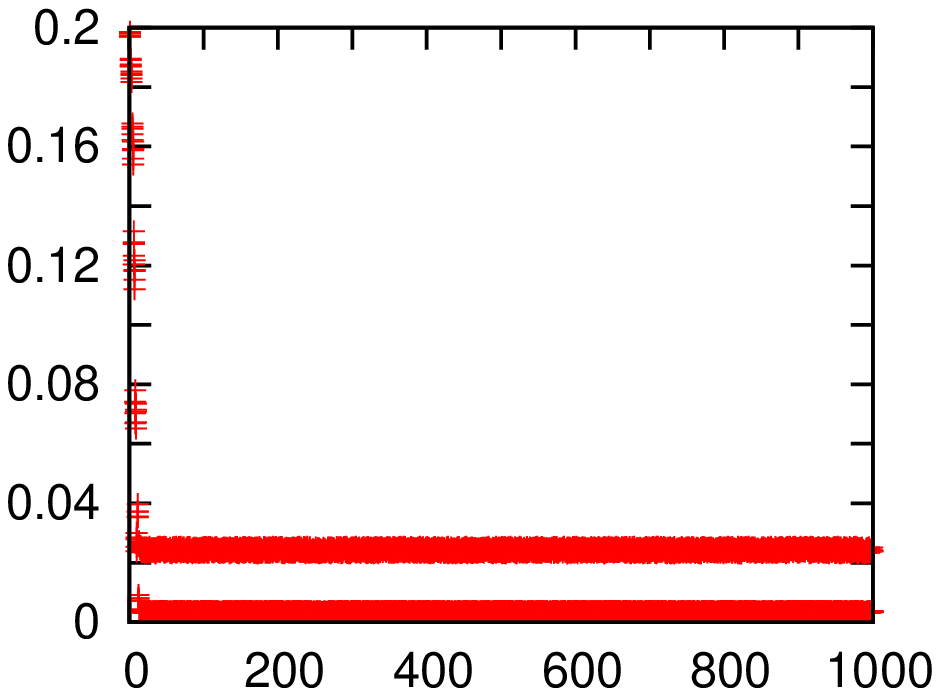}
%???????????????????????????????????
~~~%\hspace{6mm}
%!!!!!!!!!!!!!!!!!!!!!!!!!!!!!!!!!!!!!!!!!!!!
\includegraphics[height=.11\textheight]{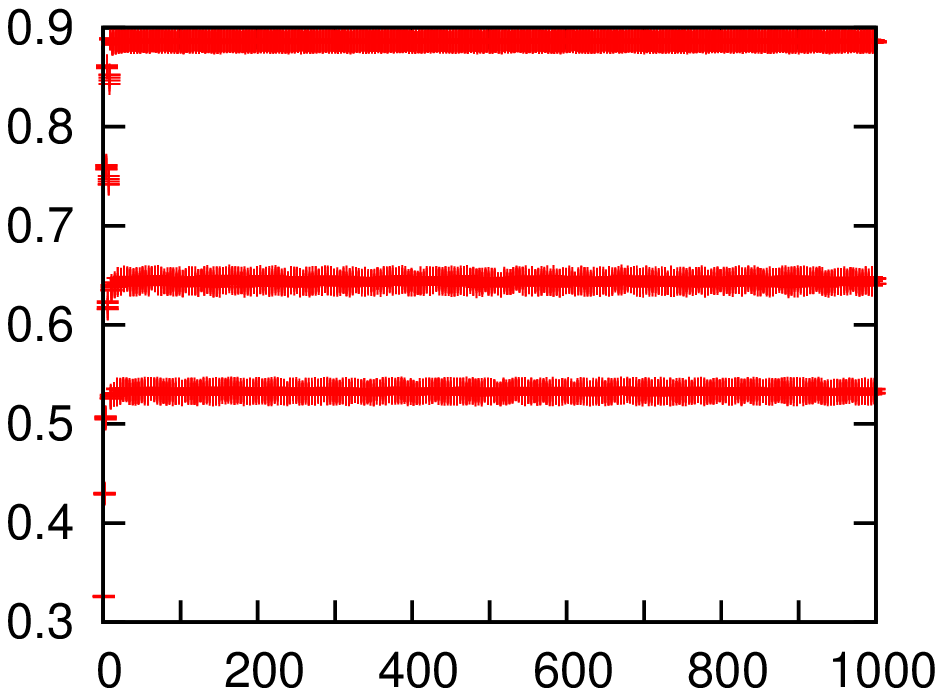}
%???????????????????????????????????
%\\
%\includegraphics[height=.12\textheight]{Prog/Impulse/Allee_6_mult_nu_0_7_cycle_2_l_1_0_01_x_0_2.ps}
%~~~%\hspace{6mm}
%\includegraphics[height=.12\textheight]{Prog/Impulse/Allee_6_mult_nu_0_7_cycle_2_l_1_0_01_x_0_3.ps}
\caption{Solutions of difference equation (\protect{\ref{eq5_7}}) with $f$ as in (\protect{\ref{def:efx6}}), PF control applied every second step,
multiplicative noise with $\ell_1=0.001$, additive noise with $\ell_2=0.01$, $\nu=0.7$ and
(left) $x_0=0.2$, (right) $x_0=0.3$.
%Lower row: five runs of the difference equation  with $f(x)= 6x^2(1-x)$, PF control applied every second step,
%multiplicative noise with $\ell_1=0.01$, $\nu=0.7$ and
%(left) $x_0=0.2$, (right) $x_0=0.3$.
}
\label{figure11}
\end{figure}

\end{example}

%%%%%%%%%%%%%%%%%%%

\section{Summary and discussion}
\label{sec:sum}

First of all, numerical simulations show less restrictive conditions on $\nu$ in \eqref{eq:TOC} than for classical (non-shifted)
PF control. If we denote $c:=1-\nu$ in \eqref{eq:TOC} then it becomes a particular case of Target Oriented Control
with an unstable equilibrium $K_1$ as a target \cite{Dattani,TPC}. 

Possible generalizations and extensions of the present research include the following topics.
\begin{enumerate}
\item [(a)]
Everywhere in simulations we assumed uniform continuous distribution, and all the estimates were dependent only on the noise
amplitude.
Specific estimates for particular types of noise distribution can be established.
\item  [(b)]
Everywhere we investigated asymptotic properties of solutions.
However, analysis of so called  transient behaviour, describing the speed of this convergence, 
starting from the initial point, maximal amplitudes for given initial values and noise characteristics, is %especially 
interesting for %real world 
applications.
%\item  [(c)]
%The present study can be extended to the case when  unbounded, for example, normal distributions are involved.
\end{enumerate}

%%%%%%%%%%%%%%

%\section{Acknowledgments} 
%
%The first author was supported by NSERC grant RGPIN-2015-05976.
%The second  and the fourth authors were supported by the Grant FEKT-S-17-4225 of Faculty of 
%Electrical Engineering and Communication, Brno University of Technology. The third  
%author was supported  by the project International Mobility of Researchers of Brno
%University of Technology  CZ.02.2.69/0.0/0.0/16-027/0008371. 

%%%%%%%%%%%%%%%%%%%%%%%%%%%

\end{document}